\documentclass[a4paper]{article}
\usepackage[latin1]{inputenc}  \usepackage[T1]{fontenc}
\usepackage{lmodern}           

\usepackage[italian,spanish,german,frenchb,english]{babel}
\usepackage{theorem}
\usepackage{amssymb}
\usepackage{amsfonts}
\usepackage{amsmath}
\usepackage{amsxtra}

\usepackage{pstricks}
\usepackage{showidx}
\usepackage{amscd}
\usepackage[active]{srcltx}
\usepackage{multicol}
\usepackage{fancyhdr}
\usepackage{changebar}
{\theoremstyle{change} \theoremheaderfont{\normalfont\bfseries}
\theorembodyfont{\slshape}
\newtheorem{Prop}{Proposition:}[section]}

{\theoremstyle{change} \theorembodyfont{\slshape}
\newtheorem{Theo}[Prop]{Theorem:}}

{\theoremstyle{change} \theorembodyfont{\slshape}
\newtheorem{Cor}[Prop]{Corollary:}}

{\theoremstyle{change} \theorembodyfont{\slshape}
\newtheorem{Lem}[Prop]{Lemma:}}

{\theoremstyle{change} \theorembodyfont{\upshape}
\newtheorem{Rem}[Prop]{\normalfont\scshape {Remark:}}}

{\theoremstyle{change} \theorembodyfont{\upshape}
\newtheorem{Defi}[Prop]{\normalfont\scshape{Definition:}}}

{\theoremstyle{change} \theorembodyfont{\upshape}
\newtheorem{Exa}[Prop]{\normalfont\scshape{Example:}}}

{\theoremstyle{change} \theorembodyfont{\upshape}
\newtheorem{Not}[Prop]{\normalfont\scshape{Notation:}}}

{\theoremstyle{change} \theorembodyfont{\upshape}
\newtheorem{Nr}[Prop]{\kern -1ex}}

{\theoremstyle{change} \theorembodyfont{\upshape}
}

\newcommand{\qed}{\hfill \mbox{\raggedright \rule{.07in}{.1in}}}

\catcode`\á=\active \def á{\'a}
 \catcode`\Á=\active \def Á{\'A}
  \catcode`\ó=\active \def ó{\'o}
  \catcode`\é=\active \def é{\'e}
  \catcode`\ú=\active \def ú{\'u}
  \catcode`\í=\active \def í{\'{\i}}
  \catcode`\ñ=\active \def ñ{\~n}
  \catcode`\Ñ=\active \def Ñ{\~N}
  \catcode`\¿=\active \def ¿{?`}
  \catcode`\º=\active \def º{$^{\underline{o}}$}
  \catcode`\ª=\active \def ª{$^{\underline{a}}$}
  \catcode`\¡=\active \def ¡{!`}
  \catcode`\â=\active \def â{\^{a}}
  \catcode`\ê=\active \def ê{\^{e}}
  \catcode`\î=\active \def î{\^{\i}}
  \catcode`\ô=\active \def ô{\^{o}}
  \catcode`\û=\active \def û{\^{u}}
  \catcode`\ç=\active \def ç{\c{c}}
  \catcode`\ü=\active \def ü{\"{u}}
  \catcode`\ö=\active \def ö{\"{o}}

 \newcommand{\ox}{\overline{x}}
 \newcommand{\yy}{\overline{y}}

 \newcommand{\ee}{\mathbb{E}}

 \newenvironment{dem}{\noindent{\it Proof}.}{\qed}

\newcommand{\nn}{\mathbb{N}}

\newcommand{\clus}{\mathcal{C}}

\newcommand{\pro}{\mathbb{P}}
\newcommand{\oo}{\mathcal{O}}

\newcommand{\mm}{\mathfrak{m}}
\newcommand{\ago}{\mathfrak{a}}

\title{Curvettes and clusters of infinitely near points}
\author{Julio Jos\'e Moyano-Fern\'andez \thanks{Research partially supported by the Deutsche Forschungsgemeinschaft (DFG), Junta de
Castilla y Le\'on grant JCyL-VA025A07 and by the Spanish
Ministerio de Educación y Ciencia grant MTM2007-64704 in the
framework of the European founds FEDER.}}

\date{Institut f\"ur Mathematik, Universit\"at
Osnabr\"uck\\
Email:  jmoyanof@uni-osnabrueck.de}
\begin{document}
\maketitle

\renewcommand{\abstractname}{Abstract}
\begin{abstract}
The aim of this paper is to revise the theory of clusters of
infinitely near points for arbitrary fields. We describe in
particular the intersection matrix of such a cluster, we introduce
the notion of curvette over an arbitrary field and we relate it to
the Hamburger-Noether tableaux associated with curves.
\end{abstract}

\noindent
AMS-Classification: 13H05, 14H20\\
Keywords: Two-dimensional regular local ring, infinitely near
point, proximity, intersection matrix, curvette

\section{Introduction}

The theory of infinitely near points was nicely introduced in the
classical treatise of Enriques and Chisini (\cite{enriques}) from
a purely geometrical point of view, based on the Max Noether's old
works. Since then, many authors have considered its algebraic
counterpart, being remarkable the works of Zariski and Lipman on
the theory of complete ideals (see \cite{zar}, \cite{lipman0},
\cite{lipman1}). Recently, these two directions have been compiled
by Casas (\cite{casas}) and Kiyek and Vicente (\cite{kiyek}).
\medskip

Infinitely near points have been mainly used for studying the
singularities of algebraic curves and their resolutions, a very
interesting subject with connections to fibre spaces, knot theory
and commutative algebra. Roughly speaking, if we blow up a closed
point $P$ on a surface $S$, we create a new surface $S^{\prime}$
containing a whole curve $E$ (called exceptional divisor) at which
$P$ used to be. Notice also that the points on $E$ are nothing but
the tangent directions at $P$ to $S$ and they are precisely the
infinitely near points to $P$. Also, it is sometimes useful to
consider normal-crossing curves at smooth points of $E$ (see for
instance \cite{de}, \cite{duke}, \cite{cadegu11} or
\cite{moyano}), the so-called curvettes---terminology introduced
by Deligne in \cite[p.13]{deligne}.
\medskip

The aim of this paper is to describe some aspects of the theory of
clusters of infinitely near points and curvettes from the
algebraic viewpoint, a topic not totally covered by the literature
as presented here. It is not our purpose to study geometric
aspects but to stress the validity of the techniques for arbitrary
fields. Nevertheless, purely algebraic objects as valuations, and
even a particular case of curvettes---the so-called general
elements---have been already studied in the investigation of the
theory of simple complete ideals of two-dimensional regular local
rings (see \cite{spiva}, \cite{grecokiyek}, \cite{grecokiyek2},
\cite{km}). We will use in this paper some of the terminology and
objects provided there, such as the Hamburger-Noether tableaux.
\medskip

The paper goes as follows. We recall in Section \ref{sec:gral} the
main concepts and results of the theory of regular local rings of
dimension two. Our main result in this section is to describe the
discrete valuation of rank $2$ defined by a regular local
two-dimensional ring $R$ and by a homogeneous prime ideal of
height $1$ of the graded ring of $R$ in terms of some
multiplicities occurring in the set of infinitely near points. In
Section \ref{sec:prox} we state the notions of cluster of
infinitely near points and proximity matrix, the latter being a
useful tool to encode the proximity relations in the cluster
introduced by Du Val in \cite{duval}. Such a matrix has to do with
the intersection relations among components of exceptional
divisors created by successive blow-ups of closed points, as we
show in Section \ref{sec:inter}; in particular, we express the
intersection matrix in terms of the cluster (Theorem
\ref{prop:intmatrix}). Section \ref{sec:HN} is devoted to describe
some numerical invariants concerning the resolution (the so-called
characteristic data) as in \cite{russell} was done for algebroid
curves, and an appropriate machinery to read them off (the
Hamburger-Noether tableau). Finally, Sections \ref{sec:cur} and
\ref{sec:approx} are devoted to show the existence and main
properties of the curvettes
 by means of the Hamburger-Noether tableau
(cf. Theorem \ref{prop:curvette}, Corollary \ref{cor:curvette},
Proposition \ref{prop:curvette2}). The Hamburger-Noether tableau
is a device that contains the most relevant data arising from the
Hamburger-Noether algorithm proposed in \cite{russell}. It is
well-known in the study of the algebroid curves (see also
\cite{campillo} in case of algebraically closed fields; more
general set-up can be founded in \cite{rybo}). In particular, we
show in our more general context that curvettes are basically the
same objects as the approximations described by Russell in
\cite{russell} for algebroid curves (see Proposition
\ref{prop:709}, Theorem \ref{prop:711}).
\medskip

An important observation for the whole paper is that the ground
field of the curve does not play any role in most of the
reasonings we do. 
\medskip

Along this paper we will denote by $\mathbb{N}$ the set of positive integer numbers, and by $\mathbb{N}_0$ the set of nonnegative integer numbers. 

\section{Generalities on two-dimensional regular local rings}
\label{sec:gral}

Along this section we will refer to the book of Kiyek and Vicente
\cite{kiyek} as a general reference. Let $R$ be a regular local
ring of dimension two with maximal ideal $\mm_R=\mm$ and residue
field $k_R$. Let $\{x,y\}$ be a regular system of para\-meters of
$R$, and let $\mathcal{K} = \mathrm{Quot} (R)$ be the field of
fractions of $R$.

\Nr For every $f \in R \setminus \{0 \}$ we define the
\textbf{order function} of $f$ as
\[
\mathrm{ord}_R (f)=\mathrm{ord} (f)=m ~ ~ ~ ~ ~ \mathrm{~if~} ~ ~
~ ~ f \in \mm^{m}, ~ f \notin \mm^{m+1}.
\]

If $m=\mathrm{ord}_R (f)$, then the class of $f$ in $\mm^m /
\mm^{m+1}$, denoted by $\mathrm{In}(f)$, is called the initial
form or the leading form of $f$. We define also the order of a
non-zero ideal $\ago$ of $R$ to be
\[
\mathrm{ord}_R (\ago)=\mathrm{ord}(\ago) := \min \{\mathrm{ord}
(a) \mid a \in \ago \}.
\]

The canonical extension of the order function to $\mathcal{K}
\setminus \{ 0\}$ gives rise to a discrete valuation of rank $1$
of $\mathcal{K}$, which we write $v_R=v$. This valuation is
non-negative on $R$ and has center $\mm$ in $R$. The valuation
ring of $v_R$ is denoted by $V_R=V$.

\Nr \label{2:2} Let $\mathcal{R}(\mm,R):=\bigoplus_{n \ge 0} \mm^n
T^n \subset R[T]$ be the Rees ring of $R$ with respect to $\mm$
for an indeterminate $T$, and let
$\mathrm{gr}_{\mm}(R):=\bigoplus_{n \ge 0} \mm^n / \mm^{n+1}$ be
the graded associated ring of $R$. Consider the homomorphism $
\varphi: \mathcal{R}(\mm,R) \to \mathrm{gr}_{\mm}(R)$. We see
immediately that $\mathrm{gr}_{\mm}(R)=k_R
[\overline{x},\overline{y}]$, where $\overline{x}:=x \mod \mm^2$
and $\overline{y} := y \mod \mm^2$, and that
$\overline{x},\overline{y}$ are algebraically independent over
$k_R$.
\medskip

Let $\pro_R$ be the set of closed points of
$\mathrm{Proj}(\mathrm{gr}_{\mm}(R))$ (i.e., homogeneous prime
ideals of $\mathrm{gr}_{\mm}(R)$ of height $1$). For $p \in
\pro_R$, the ideal $p$ is principal and generated by an
irreducible homogeneous polynomial $\overline{f} \in
k_R[\overline{x},\overline{y}]$. We set $\deg (p):= \deg
(\overline{f})$. Let $p=(\overline{f}) \in \pro_R$, where
$\overline{f} \in \mathrm{gr}_{\mm}(R)$ is homogeneous of degree
$m$, and choose $f \in \mm^m$ with $\overline{f}=f ~ \mathrm{mod~}
\mm^{m+1}$. Define $\mathfrak{n}^{\prime}_{p}:= \varphi^{-1}(p)$.
Then $\mathfrak{n}^{\prime}_{p}$ is a closed point of
$\mathrm{Proj}(\mathcal{R}(\mm,R))$ and $\mathrm{ord}_R (f)=m$.
Without loss of generality we assume that $\overline{x}$ does not
divide $\overline{f}$. Then $xT \notin \mathfrak{n}^{\prime}_{p}$
and in the ring $A:= \mathcal{R}(\mm,R)_{(xT)} =R \left [
\frac{y}{x} \right ]$ one has that the maximal ideal
$\mathfrak{n}_p$ of $A$ determined by $\mathfrak{n}_p^{\prime}$ is
$\left ( x, \frac{f}{x^{m}} \right)$. Then
$\mathcal{R}(\mm,R)_{\mathfrak{n}^{\prime}_{p}} =
A_{\mathfrak{n}_p}$ and $S_p:=A_{\mathfrak{n}_p}$ is a regular
local ring of dimension $2$ with quotient field $\mathcal{K}$ and
maximal ideal generated by $x$ and $\frac{f}{x^m}$.

\begin{Defi} \label{defn:quadrtransform}
The local ring $S_p$ is the \textbf{quadratic transform} of $R$ at
$p$. The set $N_1(R):=\{S_p \mid p \in \pro_R \}$ of all quadratic
transforms of $R$ is called the \textbf{first neighbourhood} of
$R$. Recursively, for $i > 1$, the \textbf{$i$-th neighbourhood}
of $R$, denoted by $N_{i}(R)$, is defined to be the set of
quadratic transforms of the rings in the $(i-1)$-th neighbourhood
of $R$. We also define $N(R):= \cup_{i \in \mathbb{N}_0} N_i (R)$,
i.e., the set of all two-dimensional regular local subrings of $\mathcal{K}$
containing $R$.
\end{Defi}

\Nr Let $\Omega(\mathcal{K})$ be the set of all two-dimensional
regular local subrings of $\mathcal{K}$ having $\mathcal{K}$ as
field of fractions. The elements of $\Omega(\mathcal{K})$ will be
called \emph{points}.

\Nr Let $R \in \Omega(\mathcal{K})$. If $S \in
\Omega(\mathcal{K})$ and $S \supset R$, then $S$ is said to be
\textbf{infinitely near} to $R$. In such a case there exists a
uniquely determined strictly increasing sequence
\[
R=:R_0 \subsetneqq R_1 \subsetneqq \ldots \subsetneqq R_n := S,
\eqno{(\dag)}
\]
in which $R_i \in \Omega(\mathcal{K})$ and $R_i$ is a quadratic
transform of $R_{i-1}$, for every $i \in \{1, \ldots, n \}$. In
particular, $S$ dominates $R$ and the degree extension
$[S:R]:=[k_{S}:k_{R}]$ is finite (cf. \cite{ab}; also
\cite[Chapter VII, (6.4)]{kiyek}). The previous sequence is said
to be the \textbf{quadratic sequence} between $R$ and $S$. The
integer number $n$ is called the \textbf{length} of the sequence.
Note that, if $R=S$, then we have a quadratic sequence of length
$0$.

\begin{Defi} \label{defn:proximate}
Let $S$ be an infinitely near point to the point $R$, and consider
the quadratic sequence (\dag) between $R$ and $S$. We say that $S$
is \textbf{proximate} to $R$, and we write $S \succ R$, or $R
\prec S$, if the discrete valuation ring $V_R$ contains $S$.
\end{Defi}

\Nr If $A \subset B$ are factorial rings with
$\mathrm{Quot}(A)=\mathrm{Quot}(B)$, then we associate with an
ideal $\mathfrak{a}$ of $A$ different from $0$ an ideal
$\mathfrak{a}^{B}$ in $B$, which is called the strict transform
(or ideal transform) of $\mathfrak{a}$ in $B$. For the exact
description, we refer to \cite[Chapter VII, (1.4)]{kiyek}.
\begin{Rem} \label{lem:18curvetas18}
Let $A$ be a factorial ring with quotient field $L$, and let $B
\subseteq C$ be factorial subrings of $L$ with $A \subseteq B$.
Let $\mathfrak{a}$ and $\mathfrak{b}$ be non-zero ideals in $A$.
By \cite[Chapter VII, (1.5)]{kiyek}, the following properties
hold:
\begin{enumerate}
    \item $(\mathfrak{a}^B)^C=\mathfrak{a}^C$;
    \item $(\mathfrak{a} \mathfrak{b})^B = \mathfrak{a}^B
    \mathfrak{b}^B$;
    \item if $\mathfrak{a}$ is a principal prime ideal, then either
    $\mathfrak{a} B \cap A = \mathfrak{a}$, in which case $\mathfrak{a}^B$ is a principal
    prime ideal of $B$ with $\mathfrak{a}^B \cap A = \mathfrak{a}$, or $\mathfrak{a} B \cap
    A \ne \mathfrak{a}$, in which case we have $\mathfrak{a}^B =B$.
\end{enumerate}
\end{Rem}

\Nr Let $f \in R \setminus\{0\}$. Consider the ideal $fR$ and let
$S_p$ be the quadratic transform of $R$ at $p$ (cf. \ref{2:2},
\ref{defn:quadrtransform}). Any generator of the ideal
$(fR)^{S_p}$ is called the \textbf{strict transform} of $f$ in
$S_p$. Next lemma will be needed in the sequel (cf. \cite[Chapter
VII, (2.11)]{kiyek}):

\begin{Lem} \label{lem:stricttransfisirred}
Let be the ring $R$, $p \in \pro_R$ and $S:=S_p$. Assume $\mm S =
x S$. Then we have:
\begin{itemize}
\item[(i)] If $h \in R$ is irreducible and $m:= \mathrm{ord}(h)$,
then $(hR)^S=x^{-m} h S$ and then, either $x^{-m}h$ is irreducible
in $S$ (in this case $\mathrm{In}(h) \in p $), or $x^{-m}h$ is a
unit of $S$ (and $\mathrm{In}(h) \notin p$).

\item[(ii)] Let $f,g \in R$ be irreducible and not associated (two
elements $f,g \in R$ are said to be associated if $f=ug$, where
$u$ is a unit of $R$). If $(fR)^S$, $(gR)^S$ are prime ideals of
$S$, then $(fR)^S \ne (gR)^S$ and $xS \ne (fR)^S$.
\end{itemize}
\end{Lem}

\begin{Not}
A curve $E$ in $R$ is a non-zero principal ideal $fR$ of $R$. The
element $f$ is uniquely determined up to units, and every
generator of the ideal $fR$ is called an equation of $E$. If
$fR=R$, then the curve $E$ is called empty. A curve $E$ with
equation $f$ is called irreducible, if $fR$ is a prime ideal of
$R$. Since $R$ is factorial, $fR$ is a prime ideal of $R$ if and
only if $f$ is an irreducible element of $R$. Let $E$ be a
non-empty curve with equation $f$. Let $f=f_1^{e_1} \cdot \ldots
\cdot f_r^{e_r}$ be the prime decomposition of $f$. For every $i
\in \{1, \ldots ,r \}$, let $E_i$ be the curve with equation
$f_i$. The curves $E_1, E_2, \ldots , E_r$ are called the
irreducible components of $E$, and for every $i \in \{1, \ldots ,
r \}$ $E_i$ is called irreducible component of multiplicity $e_i$.
An irreducible component of $E$ is called simple, if it has
multiplicity $1$.
\end{Not}

\begin{Prop} \label{cor:regular}
If $x \in R$ with $\mathrm{ord}_R (x)=1$, then $x$ is a regular
parameter of $R$.
\end{Prop}

\dem~ Since $x \notin \mathfrak{m}_R^2$, the element $x
\mathrm{~mod~} \mathfrak{m}_R^2$ is different from $0$ in
$\mathfrak{m}_R / \mathfrak{m}_R^2$. Therefore $x \mathrm{~mod~}
\mathfrak{m}_R^2$ takes part of a basis of the $k_R$-vector space
$\mathfrak{m}_R / \mathfrak{m}_R^2$. By \cite[Chapter IV, Korollar
2.4(b)]{kunz}, the assertion follows. \qed

\begin{Defi}
A curve $E$ in $R$ with equation $f$ is said to have no
singularities, if $\mathrm{ord}_R(f)=1$.
\end{Defi}

\begin{Rem} \label{rem:nosing}
By Proposition \ref{cor:regular}, a curve $E$ has no singularities
if and only if $f$ is a regular parameter of $R$. Consequently, a
curve with no singularities is irreducible.
\end{Rem}

\Nr \label{section:intmult} Let $f,g \in R$. We define the
intersection multiplicity between the curves $f$ and $g$ in $R$ as
the length (as $R$-modules)
\[
\iota_R (fR,gR):=\ell_R \left ( R/fR + gR \right ).
\]
This is finite if either the ideal $fR + gR$ is $\mm$-primary or
one of the elements $f,g$ is a unit in $R$ (cf. \cite[Chapter VII,
(8.6)]{kiyek}).
\medskip

Let $S \in N_1(R)$. Let $f \in R \setminus \{0\}$ be a curve such
that $(fR)^S \ne S$ and consider the exceptional divisor $\mm_R
S$. Since $\ago^S \nsubseteq \mm_R S$ for every ideal $\ago$ in
$R$, in particular the ideal $\mm_R S + (fR)^S$ is $\mm_S$-primary
and it makes sense to speak about the \emph{finite} intersection
multiplicity between the curves given by $(fR)^S$ and $\mm_R S$.
\medskip

Two curves $f,g \in R$ are said to \emph{meet transversally} at
$R$ (or at $\mm_R$) if they have no singularities and $\iota_R
(fR, gR)=1$.

\Nr The curve defined by the ideal
\[
E_{R_1}:=\mathfrak{m}_{R_0} R_1
\]
is called the exceptional divisor in $R_1$. In general, the curve
defined by the ideal
\[
E_{R_i}:=E_{R_{i-1} R_i} \cdot (E_{R_{i-1}})^{R_i}
\]
is called the \textbf{exceptional divisor} in $R_i$ for every $i
\in \mathbb{N}$ with $i > 1$, where $E_{R_{i-1}
R_i}=\mathfrak{m}_{R_{i-1}}R_i$ and $(E_{R_{i-1}})^{R_i}$ is the
strict transform of $E_{R_{i-1}}$ in $R_i$.

\begin{Lem} \label{lem:cortetransversal}
Let be the quadratic sequence between $R$ and $S$ given by
$(\dag)$. Let $i \ge 2$. Then the exceptional divisor in $R_i$
consists on either one or two components with no singularities,
and when there are two, they meet transversally at the point
corresponding to the ideal $\mm_{R_i}$.
\end{Lem}
\medskip

\dem~Let $\{x_0, y_0 \}$ be a regular system of parameters of
$R_0$. Then $\mm_{R_0}=(x_0, y_0)$ and consider the exceptional
divisor $\mm_{R_0}R_1=x_0 R_1$. Two different cases arise,
depending on whether $S \succ R$ or $S \nsucc R$.
\medskip

\emph{\textbf{Case A: $R_n=S$ is proximate to $R$.}} It means that
$R \subset S \subset V_R$ and we may choose a regular system of
parameters $\{x_1,y_1 \}$ of $R_1$ with $x_1=x_0$, $v_R (y_1)=0$
and $\mm_{R_1}=(x_1,y_1)$. Consider the quadratic transform $R_2$
of $R_1$. The exceptional divisor in $R_2$ has two components,
namely $E_{R_1 R_2} = y_1 R_2$ and $(E_{R_1})^{R_2} =
\frac{x_1}{y_1} R_2$, which meet transversally at the point
corresponding to $\mm_{R_2}= \left ( y_1, \frac{x_1}{y_1} \right
)$.
\medskip

Consider now the transforms of the exceptional divisor $E_{R_3}$
in $R_3$, i.e.: $(E_{R_1})^{R_3} = \frac{x_1}{y_1^2}R_3$,
$(E_{R_2})^{R_3} = \frac{y_1}{y_1} R_3 = R_3$ and $E_{R_2 R_3} =
y_1 R_3$. The transform $(E_{R_2})^{R_3}$ is the whole ring and
only the components $(E_{R_1})^{R_3}$ and $E_{R_2 R_3}$ survive,
and they intersect transversally at the point corresponding to the
ideal $\mm_{R_3}=\left ( y_1, \frac{x_1}{y_1^2} \right)$. We can
repeat this reasoning to show that, for $i \ge 2$, the two only
components of the exceptional divisor surviving are
$(E_{R_1})^{R_i} = \frac{x_1}{y_1^{i-1}}R_i$ and $E_{R_{i-1} R_i}
= y_1 R_i$, which meet transversally at the point given by
$\mm_{R_i}= \left ( y_1, \frac{x_1}{y_1^{i-1}} \right )$.
\medskip

\emph{\textbf{Case B:  $R_n=S$ is not proximate to $R$.}} Assume
we have the quadratic sequence
\[
R=R_0 \subset R_1 \subset \ldots \subset R_{h-1} \subset R_h
\subset R_{h+1} \subset \ldots \subset R_n,
\]
with $R \prec R_h$ and $R \nprec R_{h+1}$, for $h \ge 1$. If
$h=1$, then $\mm_{R_0}R_1 = x_1 R_1$. If $h \ge 2$, then the
exceptional divisor in $R_h$ has two components (cf. Case A);
namely $(E_{R_1})^{R_h} = \frac{x_1}{y_1^{h-1}}R_h$ and
$E_{R_{h-1} R_h} = y_1 R_h$, which intersect transversally at the
point corresponding to $\mm_{R_h}=(x_h,y_h)$, where $x_h=y_1$ and
$y_h = \frac{x_1}{y_1^{h-1}}$. We now turn to the transforms of
the exceptional divisor in $R_{h+1}$. Since $R_{h+1}$ is not
proximate to $R$, we have two possibilities:

\begin{itemize}
    \item[1)] If $E_{R_h R_{h+1}} = y_h R_{h+1}$.

    Then $\mm_{R_{h+1}}= \left ( y_h, \frac{f(x_h,y_h)}{y_h^l} \right
)$, where $f \in R, \mathrm{ord}_{R_h}(f)=l$ and $f \mathrm{~mod~}
\mm_{R_h}^{l+1} $ is an homogeneous polynomial of degree $l$. The
components of the exceptional divisor in $R_{h+1}$ are
\begin{displaymath}
\begin{array}{lclcl}
(E_{R_1})^{R_{h+1}} & = & (y_h R_h)^{R_{h+1}} & = & \left (\frac{x_1}{y_1^{h-1}}R_h \right )^{R_{h+1}} = R_{h+1} \\
(E_{R_{h}})^{R_{h+1}} & = & (x_h R_h)^{R_{h+1}} & = & \frac{x_h}{y_h} R_{h+1} \\
E_{R_h R_{h+1}} &=& y_h R_{h+1} & & .
\end{array}
\end{displaymath}
Taking into account the form of $f(x_h,y_h)$, we have

\begin{itemize}
    \item[i)] If $(f(x_h,y_h))=(x_h)$, then we have the transforms $(E_{R_{h}})^{R_{h+1}} = \frac{x_h}{y_h}R_{h+1}$ and $E_{R_{h}
R_{h+1}} = y_h R_{h+1}$, which meet transversally at the point
given by $\mm_{R_{h+1}}= \left ( y_h, \frac{x_h}{y_h} \right)$.

    \item[ii)] If $(f(x_h,y_h)) \ne (x_h)$, then $\frac{x_h}{y_h}$ is a unit in $R_{h+1}$, and only
     the component $y_h R_{h+1}$ of the exceptional divisor survives. Hence there is no intersection.
\end{itemize}

    \item[2)] If $E_{R_h R_{h+1}}= x_h R_{h+1}$, then, by the same reasoning as in the previous case B.1.), we have $\mm_{R_{h+1}}= \left ( x_h,
\frac{f(x_h,y_h)}{x_h^l} \right)$ and the components of the
exceptional divisor in $R_{h+1}$ are
\begin{displaymath}
\begin{array}{lclclcl}
(E_{R_1})^{R_{h+1}} & = & \left (\frac{x_1}{y_1^{h-1}}R_h \right )^{R_{h+1}}  & = & (y_h R_h)^{R_{h+1}} & = & \frac{y_h}{x_h} R_{h+1} \\
(E_{R_{h}})^{R_{h+1}} & = & (x_h R_h)^{R_{h+1}} & = & \frac{x_h}{x_h} R_{h+1} & = & R_{h+1} \\
E_{R_h R_{h+1}} &=& x_h R_{h+1} & & & &.
\end{array}
\end{displaymath}

We distinguish again the following two cases:

\begin{itemize}
    \item[i)] If $(f(x_h,y_h))=(y_h)$, then the components of the
    exceptional divisor are
$(E_{R_{1}})^{R_{h+1}} =\frac{y_h}{x_h}R_{h+1}$ \\
$E_{R_{h} R_{h+1}} = x_h R_{h+1}$, and they meet transversally at
a point given by the maximal ideal $\mm_{R_{h+1}}= \left ( x_h,
\frac{y_h}{x_h} \right)$.

    \item[ii)] If $(f(x_h, y_h)) \ne (y_h)$, then $(E_{R_1})^{R_{h+1}} = \frac{y_h}{x_h} R_{h+1}$; but
$\frac{y_h}{x_h}$ is a unit in $R_{h+1}$ and therefore there is no
intersection.
\end{itemize}
\end{itemize}
Notice that the components of the exceptional divisor in $R_i$
have no singularities, then they are irreducible by Remark
\ref{rem:nosing}. \qed

\begin{Prop} \label{prop:nointersection}
Let $R,S \in \Omega (\mathcal{K})$ be two points with $R \ne S$
and $R\prec S$. Let $S^{\prime} \in N_1 (S)$ with $R \prec
S^{\prime}$ and $S^{\prime \prime} \in N_1(S^{\prime})$. Then the
exceptional divisor in $S^{\prime \prime}$ consists of two
irreducible components, and they meet (transversally) if and only
if $S^{\prime \prime}$ is proximate either to $R$ or $S$.
\end{Prop}

\dem~This proposition is an easy consequence of the previous Lemma \ref{lem:cortetransversal}.  \qed

\begin{Rem}
Proposition \ref{prop:nointersection} is a generalisation for a
non-algebraically closed ground field of \cite[Proposition
4.4.2]{casas}.
\end{Rem}

\begin{Defi}
Let $n \in \mathbb{N}$, $S \in N_n(R)$, and consider the sequence
$(\dag)$ of quadratic transformations between $R$ and $S$. The
ring $S$ is said to be free with respect to $R$ if $R_{n-1}$ is
the unique ring with $S \succ R_{n-1}$; otherwise $S$ is called
satellite with respect to $R$.
\end{Defi}

From the previous facts we conclude the following result.

\begin{Cor}
Let $S \in N_n(R)$ with $n \in \mathbb{N}$ and consider the
sequence $(\dag)$ of quadratic transformations between $R$ and
$S$. Let $\{x_n,y_n\}$ be the regular system of parameters of
$R_n$ obtained from the above procedure. We have:
    \begin{enumerate}
        \item If $E_{R_n}=x_n R_n$, then the ring $R_n$ is free
        with respect to $R$; if $E_{R_n}=x_n y_n R_n$, then the
        ring $R_n$ is satellite with respect to $R$.
        \item If $E_{R_n}$ is a curve with no singularities, then
        $R_n$ is free with respect to $R$; if $E_{R_n}$ has two
        irreducible simple components, which are curves with no
        singularities, then $R_n$ is satellite with respect to
        $R$.
    \end{enumerate}
\end{Cor}

\Nr Let $p \in \mathbb{P}_R$ and let $S_p$ be the quadratic
transform of $R$ at $p$. The ideal $\mm \cdot \ago^{S_p}$ of $S_p$
is called the \textbf{reduced total transform of $\ago$ in $S_p$}
(which in \cite{kiyek} is called simply \emph{total transform}).
Let $n>1$, $S \in N_n(R)$, and consider the quadratic sequence
$(\dag)$ between $R$ and $S$. Let $\widetilde{\ago}$ be the
reduced total transform of $\ago$ in $R_{n-1}$. The ideal
$\mm_{R_{n-1}} \widetilde{\ago}^{S}$ is said to be the
\textbf{reduced total transform of $\ago$ in $S$}. We may also
describe more precisely the reduced total transform following
\cite[Chapter VII, (8.11), p. 300]{kiyek}:
\medskip

\begin{Prop} \label{prop:totaltrans}
Let $R \subset S$, with $R \ne S$ and $R,S \in
\Omega(\mathcal{K})$, then there exists a regular system of
parameters $\{x_S,y_S \}$ of $S$ such that the reduced total
transform of any ideal $\ago$ of $R$ in $S$ has the form $x_S
y_S^e \ago^S$, with $e \in \{0,1 \}$.
\end{Prop}

\dem~We will use induction on the length $n$ of the quadratic
sequence $(\dag)$ between $R$ and $S$. If $n=1$, then the ideal
$\mm \cdot \ago^{R_1}$ is the reduced total transform of $R$ in
$R_1=S$, and the result follows. Assume that the claim is true for
$n-1$, and let $R^{\prime} \in N_{n-1}(R)$ and $S$ be a quadratic
transform of $R^{\prime}$. By induction, there exists a regular
system of parameters $\{x_{R^{\prime}},y_{R^{\prime}} \}$ of
$R^{\prime}$ so that the reduced total transform of $\ago$ in
$R^{\prime}$ has the form $x_{R^{\prime}}y_{R^{\prime}}^e
\ago^{R^{\prime}}$, $e \in \{0,1\}$. Let $\mm_{R^{\prime}}$ be the
maximal ideal of $R^{\prime}$. We have to consider two cases:
\medskip

(i) If $\mm_{R^{\prime}} S = x_{R^{\prime}} S$: we set
$x_S:=x_{R^{\prime}}$. If $\frac{y_{R^{\prime}}}{x_{R^{\prime}}}$
is not a unit of $S$, then the set $\left \{ x_S,
y_S:=\frac{y_{R^{\prime}}}{x_{R^{\prime}}} \right \}$ is a regular
system of parameters of $S$, and therefore we have that
$\mm_{R^{\prime}} \cdot (x_{R^{\prime}} y_{R^{\prime}}^e
\ago^{R^{\prime}})^S=x_S y_S^e \ago^S$. If
$\frac{y_{R^{\prime}}}{x_{R^{\prime}}}$ is a unit of $S$, we
choose $y_S \in S$ such that $\{x_S,y_S \}$ is a regular system of
parameters of $S$ and we get $\mm_{R^{\prime}} \cdot
(x_{R^{\prime}} y_{R^{\prime}}^e \ago^{R^{\prime}})^S = x_S \cdot
\ago^S$.
\medskip

(ii) If $\mm_{R^{\prime}} S = y_{R^{\prime}} S$: we set
$x_S:=y_{R^{\prime}}$. If $\frac{x_{R^{\prime}}}{y_{R^{\prime}}}$
is not a unit of $S$, then the set $\left \{ x_S,
y_S:=\frac{x_{R^{\prime}}}{y_{R^{\prime}}} \right \}$ is a regular
system of parameters of $S$, and therefore we have that
$\mm_{R^{\prime}} \cdot (x_{R^{\prime}} y_{R^{\prime}}^e
\ago^{R^{\prime}})^S=x_S y_S \ago^S$. If
$\frac{x_{R^{\prime}}}{y_{R^{\prime}}}$ is a unit of $S$, we
choose $y_S \in S$ such that $\{x_S,y_S \}$ is a regular system of
parameters of $S$ and we get $\mm_{R^{\prime}} \cdot
(x_{R^{\prime}} y_{R^{\prime}}^e \ago^{R^{\prime}})^S = x_S \cdot
\ago^S$. \qed
\medskip
The result \ref{prop:totaltrans} can be extended to the following
result \ref{coro:grecKI} by Greco and Kiyek (see \cite[page
397]{grecokiyek}):

\begin{Prop} \label{coro:grecKI}
Let $R \subsetneqq S$ in $\Omega (\mathcal{K})$, and consider the
quadratic sequence between $R$ and $S$ given by $(\dag)$. We have
$\mm_S=(x,y)$ with $\mm_{R_{n-1}}S = xS$. For every non-zero ideal
$\ago$ of $R$ we have $\ago S = x^c y^d \ago^S$ where
$c:=\mathrm{ord}_{R_{n-1}} \left (\ago R_{n-1} \right )$, and if
there exists $i \in \{0, \ldots , n-2 \}$ with $S$ proximate to
$R_i$, then we have $d:= \mathrm{ord}_{R_i} \left ( \ago R_i
\right )$, and $d:=0$ otherwise.
\end{Prop}

\begin{Defi}
The ideal $\ago S$ of the previous proposition is called the
\textbf{total transform} of $\ago$ in $S$.
\end{Defi}
Notice that the reduced total transform holds the reduced
structure of the total transform.

\Nr Let $R \in \Omega (\mathcal{K})$, $p \in \mathbb{P}_R$. Set
$R_1:=S_p$ (cf. \ref{2:2}); then there exists a unique infinite
sequence
\[
R=: R_0 \subset R_1 \subset R_2 \subset \ldots \subset V_R,
\eqno(*)
\]
where $R_{i+1}$ is a quadratic transform of $R_i$ for $i > 0$.
Moreover, the union $V_p:=\bigcup_{i \ge 0} R_i$ is a valuation
ring of $K$ dominating all rings $R_i$, $i \ge 0$, it is contained
in $V_R$ and it is of the first kind with respect to $R$ (cf.
\cite[Chapter VII, (7.2)(3)]{kiyek}).
\medskip

The point $p \in \mathbb{P}_R$ defines a valuation $\nu_p:
\mathcal{K} \to \mathbb{Z}\times \mathbb{Z}$ with value group
$\mathbb{Z}\times \mathbb{Z}$ (ordered lexicographically) as
follows. Let $f \in R \setminus \{0\}$. We take the initial ideal
$(\mathrm{In}(f))$ in the ring $\mathrm{gr}_{\mm_R}(R)$ which is
homogeneous and principal, and consider its factorisation, say
\[
(\mathrm{In}(f))=\prod_{q \in \mathbb{P}_R} q^{n_q(f)}.
\]
For every $f \in R \setminus \{0\}$ we define
\[
\nu_p (f):=(\mathrm{ord}_R (f),n_p (f))
\]
and we extend to $\mathcal{K}$ in the canonical way (see
\cite[Chapter VII,~(7.5),(7.6)(4)]{kiyek}).
\medskip

In the rest of this section we see how to compute $\nu_p (f)$ in
terms of the quadratic sequence determined by $p$, i.e, we prove
the following result:

\begin{Prop} \label{prop:nuvaluation}
Let $f \in R$, $f \ne 0$ be a non-unit. Consider the factorisation
$(\mathrm{In}(f))=p_1^{n_{p_1}(f)} \cdot p_2^{n_{p_2}(f)} \cdot
\ldots \cdot p_s^{n_{p_s}(f)}$ in $\mathrm{gr}_{\mm_R}(R)$ with
$p_i \in
    \mathbb{P}_R$ and $p_i \ne p_j$ if $i \ne j$, for $1 \le i,j \le s$. Let $d_i:=
    [S_{p_i}:R]=\deg (p_i)$ for $1 \le i \le s$. Let us denote by
    $(\mathrm{qS})_j$ the quadratic sequence determined by $p_j$,
    $1 \le j \le s$. We have
\begin{enumerate}
    \item $\mathrm{ord}_R (f)= \sum_{j=1}^{s} d_j \sum_{S \in \{T \in N(S_{p_j}) \mid T \in (\mathrm{qS})_j \}} \mathrm{ord}_S ((fR)^S)$.
    \item $n_{p_j} (f)=\sum_{S \in \{T \in N(S_{p_j}) \mid T \in (\mathrm{qS})_j \}} \mathrm{ord}_S ((fR)^S)$.
\end{enumerate}
\end{Prop}

\dem~Let $p \in \mathbb{P}_R$ and let
\[
R \subset S_p=:R_1 \subset R_2 \subset R_3 \subset \ldots \subset
V_p
\]
be the quadratic sequence determined by $S_p$, where $[R_i:R]=1$
for every $i \ge 2$. Let us take $\{x=x_1,y_1\}$ a regular system
of parameters of $R_1$. One has $\mathfrak{m}_R R_1=x R_1$. Let $i
\ge 2$. We have
\[
(\mathfrak{m}_R R_1)^{R^{i}} = \left ( \frac{x}{y_1^{i-1}}\right )
R^{i}
\]
and $\mathrm{ord}_{R_i}((\mathfrak{m}_R R_1)^{R_i})=1$.
\medskip

Let us take now $R^{\prime} \in  \Omega (R_1)$ with $R^{\prime}
\nsucc R$. Set 
\[
\alpha:= \mathrm{max} \{ i \in \mathbb{N} \mid R_i
\subset R^{\prime}, \mathrm{~for~} i \ge 2\}
\]
and consider the
quadratic sequence between $R_{\alpha}$ and $R^{\prime}$
\[
R_{\alpha} \subset R_{\alpha +1} \subset R_{\alpha +2} \subset
\ldots \subset R_{\alpha + \beta}=:R^{\prime}
\]
for some $\beta \in \mathbb{N}$. Let us look at the quadratic
transform $R_{\alpha +1}$ of $R_{\alpha}$. This has the form
either $R_{\alpha +1}=R \left [ \frac{x}{y_1^{\alpha}}
\right]_{\left (y,\frac{\bullet}{y^{\star}} \right )}$, or
$R_{\alpha +1}=R \left [ \frac{y_1^{\alpha}}{x} \right]_{\left (x
, \frac{\bullet}{x^{\star}}\right )}$ (cf. \ref{2:2}). Without
loss of generality we assume that it has the first form. In such a
case $\frac{x}{y_1^{\alpha}}$ is a unit in the ring $R_{\alpha
+1}$ and therefore
\[
(\mathfrak{m}_R R_1)^{R_{\alpha +1}} = ((\mathfrak{m}_R
R_1)^{R_{\alpha +1}})^{R_{\alpha +1}} = \left (
\frac{x}{y_1^{\alpha-1}} \cdot R_{\alpha}\right )^{R_{\alpha +1}}
=R_{\alpha +1}.
\]
We have shown: for every $R^{\prime} \in  \Omega (R_1)$ with
$R^{\prime} \nsucc R$ one has $(\mathfrak{m}_R
R_1)^{R^{\prime}}=R^{\prime}$.
\medskip

From the definition of intersection multiplicity \ref{section:intmult} and  \cite[Chapter VII, (8.8)(2)]{kiyek} follows
\[
n_p (f)= \iota_{S_p} \left ( (fR)^{S_p},\mathfrak{m}_R S_p \right
) = \iota_{R_1} \left ( (fR)^{R_1}, \mathfrak{m}_R R_1 \right ).
\]

Let us assume without loss of generality that $p_1=p$. Let $\mathrm{qS}=(\mathrm{qS})_1$ be the quadratic sequence
determined by $p$. Taking into account the previous reasonings we have
\begin{align*}
n_p(f) = & \sum_{S \in N(R_1)} [S:R_1] \mathrm{ord}_S
((fR)^{R_1})^S \mathrm{ord}_S ((\mathfrak{m}_R R_1)^S) \\
= & \sum_{S \in N(R_1)} \mathrm{ord}_S
((fR)^S) \mathrm{ord}_S (\mathfrak{m}_R S) \\
= & \sum_{\substack{S \in N(R_1) \\ S \in \mathrm{qS~}}}
\mathrm{ord}_S ((fR)^S) \mathrm{ord}_S (\mathfrak{m}_R S) +
\sum_{\substack{S \in N(R_1) \\ S \notin \mathrm{qS~}}}
\mathrm{ord}_S ((fR)^S) \mathrm{ord}_S (\mathfrak{m}_R S),
\end{align*}
where $\mathrm{ord}_S ((fR)^S)=0$ and $\mathrm{ord}_S
(\mathfrak{m}_R S)=1$, when $S \notin \mathrm{qS}$. Therefore 
\[
n_p(f)=\sum_{\{S \in N(R_1) \mid S \in \mathrm{qS}\}}
\mathrm{ord}_S (fR)^S.
\]
Since every point $p_i$ determines a quadratic sequence
$(\mathrm{qS})_i$, $1 \le i \le s$, one easily deduces that
\[
n_{p_i}(f)=\sum_{\{S \in N(S_{p_i}) \mid S \in (\mathrm{qS})_i\}}
\mathrm{ord}_S ((fR)^S) \eqno{(\ddag)}
\]
for $i \in \{1, \ldots s\}$, and (ii) is proven. The first
statement follows straightforward from ($\ddag$) and the fact that
$\mathrm{ord}_R (f)=\sum_{j=1}^{s} d_j \cdot n_{p_j}(f)$. \qed

\section{Proximity matrices for clusters of infinitely near
points}\label{sec:prox}

\Nr Let $\Omega (R)$ be the set of all two-dimensional regular
local subrings of $\mathcal{K}$ containing the ring $R$. Notice that if $S
\in \Omega (R)$, then $\mm_S \cap R = \mm_R$. Since the set $\Omega
(R)$ consists of infinitely many elements, it would be more appropriate to
deal with suitable finite subsets, which are described in the following definition.

\begin{Defi}
A \textbf{cluster} in $\Omega (R)$ over $R$, denoted by $\clus
(R)$ (or simply $\mathcal{C}$, if there is no risk of confusion),
is a finite subset of $\Omega (R)$ such that
\begin{itemize}
    \item[(i)] the point $R \in \mathcal{C}$;
    \item[(ii)] if $R^{\prime} \in \clus$ and $R=:R_0 \subset R_1 \subset \ldots \subset
    R_n:=R^{\prime}$ is the quadratic sequence between $R$ and
    $R^{\prime}$, then $R_i \in \clus$ for all $i \in \{1, \ldots, n-1 \}$.
\end{itemize}
\end{Defi}

\Nr \label{defn:proximity} Let $\mathcal{C}$ be a cluster in
$\Omega (R)$ over $R$. The \textbf{proximity matrix} associated
with $\mathcal{C}$ is defined to be the matrix
$P_{\mathcal{C}}=(p_{S,T})$, for every $S, T \in \mathcal{C}$, where
\begin{displaymath}
p_{S,T}=\left\{%
\begin{array}{rl}
    1, & \hbox{if \ $S=T$;} \\
    -1, & \hbox{if \ $S \prec T$;} \\
    0, & \hbox{otherwise.}
\end{array}%
\right.
\end{displaymath}

Consider also the diagonal matrix $\Delta_{\mathcal{C}} =
(d_{S,T})$, for every $S,T \in \mathcal{C}$, given by
\begin{displaymath}
d_{S,T}=\left\{%
\begin{array}{rl}
    [S:R], & \hbox{if \ $S=T$;} \\
    0, & \hbox{otherwise.}
\end{array}%
\right.
\end{displaymath}

The proximity matrix can be slightly turned out to a matrix
$P^{\prime}_{\mathcal{C}}:=\Delta_{\mathcal{C}}^{-1} \cdot
P_{\mathcal{C}} \cdot \Delta_{\mathcal{C}}$ with entries
$(p^{\prime}_{S,T})$, $S,T \in \mathcal{C}$, where
\begin{displaymath}
p^{\prime}_{S,T}=\left\{%
\begin{array}{rl}
    1, & \hbox{if \ $S=T$;} \\
    -[S:T], & \hbox{if \ $S \prec T$;} \\
    0, & \hbox{otherwise.}
\end{array}%
\right.
\end{displaymath}
Such a matrix was proposed by Lipman in \cite{lipman} in order to
encode the pro\-ximity inequalities in a shorter way, and it is
called the \textbf{refined proximity matrix} associated with
$\mathcal{C}$. Nevertheless, this matrix does not take into
account all possible field extensions from the origin on. To
obtain that, we introduce a matrix $\widetilde{P}_{\mathcal{C}}$
with entries $(\widetilde{p}_{S,T})$, $S, T \in \mathcal{C}$,
where
\begin{displaymath}
\widetilde{p}_{S,T}=\left\{%
\begin{array}{rl}
    [S:R], & \hbox{if \ $S=T$;} \\
    -[S:R], & \hbox{if \ $S \prec T$;} \\
    0, & \hbox{otherwise.}
\end{array}%
\right.
\end{displaymath}
We will call it the \textbf{total proximity matrix} associated
with $\mathcal{C}$.

\begin{Rem}
From Definition \ref{defn:proximate}, it is easy to check that
both the matrix $P_{\mathcal{C}}, P_{\mathcal{C}}^{\prime}$ and
$\widetilde{P}_{\mathcal{C}}$ are invertible, and the entries of
$P_{\mathcal{C}}^{-1}, (P^{\prime}_{\mathcal{C}})^{-1}$ and
$\widetilde{P}_{\mathcal{C}}^{-1}$ are non-negative integers (it
also follows from \cite[Corollary 4.6]{lipman}).
\end{Rem}

\begin{Exa} \label{exa:proxi}
Let $\mathbb{F}_2[x,y]$ be the polynomial ring of two
indeterminates over the field $\mathbb{F}_2$. Let us consider the
maximal ideal $\mm=(x,y)$ and the localisation $R_0:=R=
\mathbb{F}_2[x,y]_{\mm}$. The residue field of $R$ is
$k_R=\mathbb{F}_2$. Take a point $p \in \pro_R$ given by $\overline{f}=\overline{y}^{2}-\overline{x}^5 \in
\mathbb{F}_2[\overline{x},\overline{y}]$. Since $\overline{x}$ does not lie in $p$, the quadratic transform
of $R$ at $p$ is the local ring
\[
R_1:=S_p = R \left [ \frac{\mm}{x} \right ]_{\mathfrak{n}_p} = R
\left [  \frac{y}{x} \right ]_{\mathfrak{n}_p},
\]
where $\mathfrak{n}_p := \left ( x, \frac{y^2-x^5}{x^2} \right )$
is an ideal of $R \left [ \frac{y}{x} \right ]$ which is prime of
height $2$. The residue field of $S_p$ is $k_{S_p}=\mathbb{F}_4$.
The strict transform of $\overline{f}R$ in $S_p$ is
$\overline{f}^{(1)}=\overline{y}^2-\overline{x}^3$. The
exceptional divisor has the form $xR$ in $S_p$. Now take the point
$p^{\prime}$ given by $\overline{f}^{(1)}$. If we write again
$\mathfrak{m}^{(1)}=(x,y)$ the maximal ideal of the ring
$R^{(1)}=\mathbb{F}_4[x,y]_{\mathfrak{m^{(1)}}}$, we obtain that
$R_2:=S_{p^{\prime}}=R^{(1)}[\frac{y}{x}]_{\mathfrak{n}_{p^{\prime}}}$
with $\mathfrak{n}_{p^{\prime}} := \left ( x, \frac{y^2-x^3}{x^2}
\right )$ and $k_{S_{p^{\prime}}}$ is an extension of degree $2$
of $\mathbb{F}_4$, also it is $\mathbb{F}_6$. Again we take the
point $p^{\prime \prime}$ given by the strict transform
$\overline{f}^{(2)}=\overline{y}^2-\overline{x}$, the maximal
ideal $\mathfrak{m}^{(2)}=(x,y)$ and
$R^{(2)}=\mathbb{F}_6[x,y]_{\mathfrak{m}^{(2)}}$. Since
$\overline{y}$ does not lie in $p^{\prime \prime}$, the transform
is in this case $R_3:=S_{p^{\prime
\prime}}=R^{(2)}[\frac{y}{x}]_{\mathfrak{n}_{p^{\prime \prime}}}$
with $\mathfrak{n}_{p^{\prime \prime}} := \left ( y,
\frac{y^2-x}{y} \right )$ and
$k_{S_{p^{\prime}}}=k_{S_{p^{\prime}}}$. Now we pick the point
$p^{\prime \prime \prime}$ given by the strict transform
$\overline{f}^{(3)}=\overline{y}-\overline{x}$. Again
$\overline{y}$ does not lie in $p^{\prime \prime \prime}$ and
$R_4:=S_{p^{\prime \prime
\prime}}=R^{(3)}[\frac{x}{y}]_{\mathfrak{n}_{p^{\prime \prime
\prime}}}$ with $\mathfrak{n}_{p^{\prime \prime \prime}} := \left
( y, \frac{y-x}{y} \right )$ and no extension of the residue
field. Finally, the strict transform
$\overline{f}^{(4)}=\overline{y}-1$ is a unit and we finish. We
obtain the quadratic sequence $R=:R_0 \subseteq R_1 \subseteq R_2 \subseteq R_3 \subseteq R_4$
with the following proximity relations: $R_1 \succ R_0$, $R_2
\succ R_1$, $R_3 \succ R_1$, $R_3 \succ R_2$, $R_4 \succ R_3$. It
is easy to see that the proximity matrices corresponding to the
definitions in \ref{defn:proximity} are
\[
P=\left(%
\begin{array}{ccccc}
  1 & -1 & 0 & 0 & 0 \\
  0 & 1 & -1 & -1 & 0 \\
  0 & 0 & 1 & -1 & 0 \\
  0 & 0 & 0 & 1 & -1 \\
  0 & 0 & 0 & 0 & 1 \\
\end{array}%
\right),
P^{\prime}=\left(%
\begin{array}{ccccc}
  1 & -2 & 0 & 0 & 0 \\
  0 & 1 & -2 & -2 & 0 \\
  0 & 0 & 1 & -1 & 0 \\
  0 & 0 & 0 & 1 & -1 \\
  0 & 0 & 0 & 0 & 1 \\
\end{array}%
\right)
\]
\[
\widetilde{P}= \left(%
\begin{array}{ccccc}
  1 & -2 & 0 & 0 & 0 \\
  0 & 2 & -4 & -4 & 0 \\
  0 & 0 & 4 & -4 & 0 \\
  0 & 0 & 0 & 4 & -4 \\
  0 & 0 & 0 & 0 & 4 \\
\end{array}%
\right).
\]
\end{Exa}

\section{Intersection matrix in terms of a cluster}
\label{sec:inter}

\Nr Let $k$ be a field. Let $X_0$ be a two-dimensional regular
scheme of finite type over $k$. Take a closed point $x_0 \in X_0$
and blow up at $x_0$ to obtain another two-dimensional regular
scheme $X_1$ and repeat the process $s$ times. We get a finite
sequence of blowing ups
\[
X=X_{s} \overset{\pi_{s}}{\longrightarrow} X_{s-1}
\overset{\pi_{s-1}}{\longrightarrow} \ldots \longrightarrow X_3
\overset{\pi_3}{\longrightarrow} X_2
\overset{\pi_2}{\longrightarrow} X_1
\overset{\pi_1}{\longrightarrow} X_0, \eqno(\ddag)
\]
obtained in this way, where $\pi_i$ is the blowing up at a closed
point $x_{i-1} \in X_{i-1}$, for $1 \le i \le s$. Every point
$x_i$ provides a two-dimensional regular local ring $R_i :=
\oo_{X_{i},x_i}$, where $R_0=:R$. Notice also that $k=k_R$.
Sometimes we will speak about points $x_i$ instead of rings $R_i$
and we will apply the notations used for the rings to the points.
In particular, if a ring $R_i=\oo_{X_i,x_i}$ is proximate to a
ring $R_j=\oo_{X_j,x_j}$ for some $i,j \in \mathbb{N}_0$, then we
will write either $R_i \succ R_j$ (as in Definition
\ref{defn:proximity}), or $x_i \succ x_j$. Moreover, the ring
homomorphism $R_0:=R \longrightarrow R_i$ induces the field
extension $k_R \hookrightarrow R_i / \mm_{R_i}=:k_{R_i}$. For
convenience, we will denote the residue field $k_{R_i}$ simply by
$k_i$. The degree of this field extension is finite, and it will
be denoted by $h_i$ or $[R_i :R]$, as we have already seen.
Finally, we will set $\pi := \pi_{s} \circ \pi_{s-1} \circ \ldots
\circ \pi_2 \circ \pi_1$.

\Nr \label{defn:n} For every point $x \in X$, we have a local ring
$\oo_{X,x} \in \Omega_s (R):=N_0(R) \cup \ldots \cup N_s (R)$
which determines a quadratic sequence of length $n \le s$ between
$R$ and $\oo_{X,x}$. Furthermore, we associate to the sequence
($\ddag$) a cluster over $R$, which will be denoted by
$\mathcal{C}_s$ or $\mathcal{C}$, if no risk of confusion arises.
The diagonal matrix having the values $h_i$ on the diagonal will
be denoted by $\Delta_{\mathcal{C}}$.

\begin{Not}
We will write $E_{i,i}$ for the exceptional divisor of $\pi_i$ as
divisor of $X_i$, and we denote by $E_{i,j}$ (resp.
${E_{i,j}^{\ast}}$) the strict transform (resp. the total
transform) of $E_{i,i}$ in $X_j$ by the morphism $X_j
\longrightarrow X_i$, for $j>i$. We denote by $E_i$ (resp.
$E_i^{\ast}$) the strict (resp. total) transform $E_{i,s}$ (resp.
${E_{i,s}^{\ast}}$) by the morphism $X \longrightarrow X_i$.
\end{Not}

Let $\ee$ be the subgroup of $1$-cycles of $X$ of the form
$\sum_{i=1}^{s} n_i E_i$, with $n_i \in \mathbb{Z}$ (i.e., the
free $\mathbb{Z}$-module generated by the divisors $E_i$). Both
 $E=(E_1, \ldots , E_{s})$ and $E^{\ast}=(E^{\ast}_1, \ldots ,
 E^{\ast}_{s})$ are $\mathbb{Z}$-basis of $\ee$. More precisely,
 the proximity matrix of a cluster $\mathcal{C}$ over
 $R=\oo_{X_0,x_0}$ is the matrix of the change of basis from
$E$ to $E^{\ast}$:

\begin{Lem} \label{lemma:18}
Let $P_{\mathcal{C}}$ be the proximity matrix of a cluster
$\mathcal{C}$ over $R$. Then $E=E^{\ast} P_{\mathcal{C}}$.
\end{Lem}

\dem~We have to prove that $E_i = \sum_j p_{i,j} E_j ^{\ast}$,
where $p_{i,j}=1$ if $i=j$, $p_{i,j}=-1$ if $x_j \succ x_i$ and
$0$ otherwise. The point $x_j$ is proximate to $x_i$ if and only
if $x_j \in E_{i,j}$; in other words, $E_{i,j+1} \cap E_{j, j+1}
\ne \emptyset$. We deduce that the multiplicity of $E_{i,j}$ at
$x_j$ is $1$ and also that
\begin{align*}
 \pi_j^{\ast} E_{i,j} & =  E_{i,j+1} + E_{j,j+1} \mathrm{~if~} x_j \succ x_i \\
  \pi_j^{\ast} E_{i,j}& =  E_{i,j+1} \mathrm{~otherwise.~}
\end{align*} These two equalities combined yield
\begin{align*}
  \pi_j^{\ast} E_{i,j}& = E_{i,j+1} - p_{i,j} E_{j, j+1}.
\end{align*}

Write $\widehat{E}_{i,j}^{\ast}$ for the total transform
$(\pi^{\ast}_{j-1} \circ \ldots \circ \pi^{\ast}_{i+1})(E_{i,i})$
of $E_{i,i}$ in $X_j$. We will show by induction on $j-i$ that,
for $i < j$, the following equality holds:
\[
\widehat{E}_{i,j}^{\ast} + p_{i,i+1} \widehat{E}_{i+1,j}^{\ast} +
\ldots + p_{i,j-1} \widehat{E}_{j-1,j}^{\ast} = E_{i,j}.
\]

For $j=i+1$ it is obvious. If it holds for particular values of
$j$ and $i$, then, if we apply $\pi_j^{\ast}$ to both sides we get
\[
\widehat{E}_{i,j+1}^{\ast} + p_{i,i+1}
\widehat{E}_{i+1,j+1}^{\ast} + \ldots + p_{i,j-1}
\widehat{E}_{j-1,j+1}^{\ast} = \pi^{\ast}_{j} E_{i,j}.
\]
Hence the result follows for $i$ and $j+1$, because
\[
p_{i,j} \widehat{E}^{\ast}_{j,j+1} + \pi^{\ast}_{j} E_{i,j} =
E_{i,j+1}.
\]
Taking $j=s$ in the equation just obtained, one has
\[
E_i = E_i^{\ast} + p_{i,i+1} E_{i+1}^{\ast} + \ldots + p_{i,k}
E_k^{\ast}.
\]
Since $p_{i,i}=1$ and $p_{i,j}=0$ for $i > j$, we are done. \qed
\medskip

\Nr Moreover, on every $X_i$ occurring in the sequence of
blowing-ups ($\ddag$) we define the intersection of cycles. We
have a symmetric bilinear intersection form given as follows:
\begin{displaymath}
\begin{array}{ccc}
  \mathbb{E} \times \mathbb{E} & \longrightarrow & \mathbb{Z} \\
 (A,B) & \mapsto & (A \cdot B),
\end{array}
\end{displaymath}
i.e. it is given by intersecting cycles (cf. \cite[\S 9.1.2 and
Proposition 2.5]{liu}). If we denote by $h_i$ the degree of the
extension $k_R \subset k_i$, by the projection's formula (see
\cite[Theorem 9.2.12, p. 398]{liu}) we get
\begin{equation} \label{eq:deltakronecker}
(E_i^{\ast} \cdot E_j^{\ast}) = -\delta_{ij} h_i, \nonumber
\end{equation}
where $\delta_{ij}$ is the Kronecker's delta. Therefore the matrix
of the intersection form in the basis $E^{\ast}$ is
$-\Delta_{\mathcal{C}}$. By Lemma \ref{lemma:18}, the matrix of
the intersection form in the basis $E$ is
\[
N_{\mathcal{C}}:= -P_{\mathcal{C}} \cdot \Delta_{\mathcal{C}}
\cdot P_{\mathcal{C}}^t,
\]
with $P_{\mathcal{C}}^t$ the transpose of $P_{\mathcal{C}}$. Note
that $P_{\mathcal{C}} \cdot \Delta_{\mathcal{C}} =
\widetilde{P}_{\mathcal{C}}$, i.e. the total proximity matrix.

\begin{Defi}
The matrix $N_{\mathcal{C}}$ is called the \textbf{intersection
matrix} (with respect to the basis $E$) associated with the
cluster $\mathcal{C}$.
\end{Defi}
In other words, the intersection matrix associated with the
cluster $\clus$ is $ N_{\mathcal{C}}:=
-\widetilde{P}_{\mathcal{C}} \cdot P^{t}_{\mathcal{C}}$.

\begin{Exa}
If we take the quadratic sequence of the Example \ref{exa:proxi},
a simple calculation shows that the intersection matrix is
\begin{displaymath}
N=\left(%
\begin{array}{ccccc}
  -3 & 2 & 0 & 0 & 0 \\
  2 & -10 & 0 & 4 & 0 \\
  0 & 0 & -8 & 4 & 0 \\
  0 & 4 & 4 & -8 & 4 \\
  0 & 0 & 0 & 4 & -4 \\
\end{array}%
\right).
\end{displaymath}
\end{Exa}

We can characterise the entries of the intersection matrix of a
cluster as follows:

\begin{Theo} \label{prop:intmatrix}
The entries $n_{S,T}$ of the intersection matrix
$N_{\clus}=(n_{S,T})$, for every $(S,T) \in \clus \times \clus$ are
\begin{itemize}

    \item[*] $-\left ( [S:R] + \sum_{\substack{U \in \clus \\ U \succ S}} [U:R]\right
    )$ if  $S=T$;

    \item[*] $[T:R]$ if $T \succ S$ and the point $T^{\ast} \in
    N_1(T)$ with $S \prec T^{\ast}$ does not belong to the cluster
    $\clus$;

    \item[*] $[S:R]$ if $S \succ T$ and the point $S^{\ast} \in
    N_1(S)$ with $T \prec S^{\ast}$ does not belong to the cluster $\clus$;

    \item[*] $0$ otherwise.
\end{itemize}
\end{Theo}

\dem ~ From Definition \ref{defn:proximity}, the entries of the
matrix $N_{\clus}$ are
\[
n_{S,T} = - \sum_{U \in \clus} \widetilde{p}_{U,S} p_{U,T}.
\]
If $S=T$, then we have
\begin{align*}
  n_{S,S} & =  - \sum_{\{U \in \clus \mid U \supset S \}} \widetilde{p}_{U,S} p_{U,S} =  - \sum_{\{U \in \clus \mid U \succ S\}} \widetilde{p}_{U,S} p_{U,S} -\widetilde{p}_{S,S}p_{S,S} \\
   & =  - \sum_{\{U \in \clus \mid U \succ S\}} [U:R] - [S:R].
\end{align*}
If $S \ne T$, then we have three possibilities, namely:
\begin{itemize}
    \item[1)] $S \nsubseteq T$ and $T \nsubseteq S$.
    \medskip
    \begin{align*}
  n_{S,T} & = - \sum_{\{U \in \clus \mid U \supset S, U \supset T\}} \widetilde{p}_{U,S} p_{U,T} = 0.
\end{align*}

    \item[2)] $S \subset T$.
    \medskip
    \begin{itemize}
        \item[a)] If $S \nprec T$, then we have
\begin{align*}
  n_{S,T} & = - \sum_{\{U \in \clus \mid U \supset S, U \supset T\}} \widetilde{p}_{U,S} p_{U,T}= - \sum_{\{U \in \clus \mid U \succ T\}} \widetilde{p}_{U,S}p_{U,T}-\widetilde{p}_{T,S}p_{T,T} \\
   & =  -\sum_{\{U \in \clus \mid U \succ T\}} \widetilde{p}_{U,S}p_{U,T}-0= 0.
\end{align*}
Namely, since $U \succ T$, we have $U \subseteq V_T$; if we assume
that $U \supset S$, then $U \subseteq V_S$. Since $S \nprec T$,
$V_T \nsubseteq V_S$ and so $U \subseteq V_T \nsubseteq V_S$, then
$U \nsucc S$ and therefore $\widetilde{p}_{U,S}=0$ and
$n_{S,T}=0$.
\medskip

        \item[b)] If $S \prec T$, then there exists a
point $T^{\ast} \in N_1 (T)$ satisfying $S \prec T^{\ast}$. We
distinguish two cases:
\medskip

\begin{itemize}
    \item[i)] if $T^{\ast} \in \clus$, then
    \begin{align*}
  n_{S,T} & =  - \sum_{\{U \in \clus \mid U \supset S, U \supset T\}} \widetilde{p}_{U,S} p_{U,T}=  - \sum_{\{U \in \clus \mid U \succ T\}} \widetilde{p}_{U,S}p_{U,T}-\widetilde{p}_{T,S}p_{T,T} \\
   & =  -\sum_{\{U \in \clus \mid U \succ T\}} \widetilde{p}_{U,S}p_{U,T}+[T:R].
\end{align*}
To compute $\sum_{\substack{U \in \clus \\ U \succ T}}
\widetilde{p}_{U,S}p_{U,T}$, let us consider the quadratic
sequence
\[
R_0=R \subset \ldots \subset R_s=S \subset \ldots \subset R_t=T
\subset R_{t+1}=T^{\ast} \subset \ldots
\]
If $U=T^{\ast}$, then $U \succ T$ and $U \succ S$; therefore $U$
cannot be proximate to any other point of the sequence. Then
$\widetilde{p}_{U,S}=\widetilde{p}_{T^{\ast},S}=-[U:R]=-[T^{\ast}:R]$
and $p_{U,T}=p_{T^{\ast},T}=-1$. Whenever $U \ne T^{\ast}$,
suppose $U=R_{t+i}$ for some $i \ge 2$; then $U$ is proximate to
$R_{t+i-1}$ and proximate to $R_t=T$ as well, hence $U$ cannot be
proximate to $S$ and so $\widetilde{p}_{U,S}=0$. Then
\begin{align*}
  n_{S,T} & =  -\sum_{\{U \in \clus \mid U \succ T\}} \widetilde{p}_{U,S}p_{U,T} +[T:R] =  -\widetilde{p}_{T^{\ast},S} p_{T^{\ast},T} + [T:R]   \\
   & =  -((-[T^{\ast}:R])(-1)) + [T:R] =  -([T^{\ast}:T][T:R]) + [T:R] \\
   & =  [T:R] (1-[T^{\ast}:T]).
\end{align*}
But, by \cite[Chapter VII, (7.2)(2)]{kiyek}, we have
$[T^{\ast}:T]=1$ and therefore $n_{S,T}=[T:R](1-1)=0$.

    \item[ii)] if $T^{\ast} \notin \clus$, then consider a
    quadratic sequence as above
    \[
    R_0=R \subset \ldots \subset R_s=S \subset \ldots \subset R_t=T \subset R_{t+1}=T^{\ast} \subset \ldots
    \]
    Now we have
    \begin{align*}
  n_{S,T} & =  -\sum_{\{U \in \clus \mid U \succ T\}} \widetilde{p}_{U,S}p_{U,T}-\widetilde{p}_{T,S}p_{T,T}\\
  & = -\sum_{\{U \in \clus \mid U \succ T\}}\widetilde{p}_{U,S}p_{U,T}+[T:R].
\end{align*}

    In this case, the ring $U$ cannot be equal to $T^{\ast}$, because $T^{\ast} \notin
    \clus$ and therefore, by the same reasoning as in the case $U \ne
    T^{\ast}$ within the previous item i), we have that, for all $U \in \clus$ with
    $U \succ T$, then $\widetilde{p}_{U,S}=0$ and $n_{S,T}=[T:R]$.
\end{itemize}
\medskip

    \end{itemize}

    \item[3)] $T \subset S$.
    \medskip
    This situation is totally analogous to the previous case $S \subset T$. \qed
\end{itemize}

\begin{Rem}
From the previous arguments it is now easy to see that the
intersection matrix shows whether the components $E_i$ of the
exceptional divisors occurring in a blowing-up process intersect.
Indeed, the entries $n_{i,j}$ of the intersection matrix
$N_{\mathcal{C}_s}=(n_{i,j})$, for all $1 \le i,j \le s$, are
\[
n_{i,j} = \left \{
\begin{array}{ll}
-h_i - \sum_{p_l \succ p_i} h_l, & \mathrm{if~} i=j; \\
{}[k_P:k_R], & \mathrm{if~} i \ne j  \mathrm{~and~} E_i \cap E_j = \{P \}; \\
0, & \mathrm{if~} i \ne j  \mathrm{~and~} E_i \cap E_j = \emptyset, \\
\end{array}%
\right.
\]
where $k_P$ is the residue field of the local ring of $X$ at $P$.
\end{Rem}

\section{Hamburger-Noether Tableau and characteristic data}
\label{sec:HN}

Let $f \in R$ be an analytically irreducible curve. We want to
define the Hamburger-Noether tableau of $f$ following
\cite{russell}, \cite{grecokiyek}. The only difference with
\cite{russell} here is the use of an arbitrary field instead of an
algebraically closed one, which makes us to consider some Galois
groups of the corresponding field extensions. This will cause a
slight modification in the algorithm in \cite{russell}, as we now
briefly explain.

\Nr \label{sec:51} Let $V=k[\![t]\!]$ the ring of formal power
series over a field $k$ in the indeterminate $t$, $\mathfrak{n}$
its maximal ideal and $v$ the discrete valuation associated to
$k(\!(t)\!)$.
\medskip

Let $x,y \in \mathfrak{n}$, with $(x,y) \ne (0,0)$. We define a
matrix
\begin{equation}
\mathrm{HN}(x,y)=\left(%
\begin{array}{c}
  p_i \\
  c_i \\
  a_i \\
\end{array}%
\right)_{1 \le i < \infty} \nonumber
\end{equation}
with $p_i,c_i \in \nn \cup \{\infty\}$, $a_i \in \overline{k}
\setminus \{ 0\}$ for every $i \in \mathbb{N}$, by means of the
following algorithm (cf. \cite{russell}).
\medskip

If $x=0$, then $y \ne 0$ (since $x$ and $y$ cannot vanish
simultaneously) and set $p_i:=v(y)$, $c_i:=v(x)=\infty$ and
$a_i:=0$ for every $i \in \mathbb{N}$. If $x \ne 0$, then we put
$x_0:=x, y_0:=y$, and we set $x_1:=x, y_1=1,z_1=y$. If $z_1 \ne
0$, then we put $\eta_0 :=  y_0$ and $\eta_1 := x_1$, and we
define $\kappa \in \nn$, non-zero elements $\eta_2, \ldots ,
\eta_{\kappa+1} \in V$ and $s_1, \ldots , s_{\kappa} \in \nn_0$ by
the requirement that
\begin{eqnarray}
\eta_{i-1} = \eta_i^{s_i} \eta_{i+1} & \mathrm{~for~every~} & i \in \{1, \ldots , \kappa \} \nonumber \\
0 < v(\eta_i) < v(\eta_{i-1}) & \mathrm{~for~every~} & i \in \{2,
\ldots , \kappa \} ~ ~ \mathrm{~and~} v(\eta_{\kappa+1}) = 0.
\nonumber
\end{eqnarray}

Notice that
\[
v(\eta_{i-1})=s_i \cdot v(\eta_i) + v(\eta_{i+1}) ~ ~ ~
\mathrm{~for~} i \in \{ 1, \ldots , \kappa \}
\]
is the Euclidean algorithm for the natural integers $v(\eta_0),
v(\eta_1)$, therefore $v(\eta_{\kappa}) = \mathrm{gcd}
(v(\eta_0),v(\eta_1))$. From
$v(\eta_{\kappa-1})=v(\eta_{\kappa}^{s_{\kappa}})< \infty$ we see
that $\eta_{\kappa-1}/ \eta_{\kappa}^{s_{\kappa}}$ is a unit in
the integral closure $\overline{R}$ of $R$; then there exists a
unique $a:=a(\eta_0, \eta_1)$ which is a non-zero element of the
extension field $k_1$ of $k:=k_0$ of degree $d_1=\sharp \left (
\mathrm{Gal}(\overline{k_0}/k_0) \right )$, where
$\mathrm{Gal}(\overline{k_0}/k_0)$ denotes the Galois group of the
extension $\overline{k_0}/k_0$, such that
\[
v \left ( \eta_{\kappa-1} - a \eta_{\kappa}^{s_{\kappa}} \right )
> v (\eta_{\kappa-1}).
\]
Hence we set $a_1:=a(\eta_0,\eta_1), p_1:=v(z_1)-v(y_1),
c_1:=v(x_1)$, and also $x_2:=\eta_{\kappa}, y_2:=\eta_{\kappa-1},
z_2:=\eta_{\kappa-1}-a_1 \eta_{\kappa}^{s_{\kappa}}$. On the other
hand, if $z_1=0$ the we define $a_1:=0, p_1:=\infty, c_1:=v(x_1) <
\infty$, and $x_2:=x_1, y_2:=y_1, z_2:=z_1$. Then we apply the
algorithm again to $x_2,y_2,z_2$. In general, if we assume that
$x_i,y_i,z_i$ have been already computed for $i \in \mathbb{N}$,
then we have:

\begin{itemize}
    \item If $z_i \ne 0$, then we set $\eta_0:=z_i, \eta_1:=x_i$
    and $a_i:=a(\eta_0,\eta_1) \in k_i \setminus \{0\}, p_i:=v(z_i)-v(y_i), c_i:=v(x_i)$
    and $x_{i+1}:=\eta_{\kappa}, y_{i+1}:=\eta_{\kappa-1}, z_{i+1}:=\eta_{\kappa-1}-a_i
    \eta_{\kappa}^{s_{\kappa}}$.
    \item If $z_i=0$, then we put $a_i:=0, p_i:=\infty, c_i:=v(x_i)$
    and also $x_{i+1}:=x_i, y_{i+1}:=y_i, z_{i+1}:=z_i$.
\end{itemize}
\medskip

\begin{Defi}
Let $x,y \in \mathfrak{n}$ with $(x,y) \ne (0,0)$. The matrix
$\mathrm{HN}(x,y)$ defined in \ref{sec:51} will be called the
\emph{Hamburger-Noether tableau} of the pair $(x,y)$ in $V$ (cf.
\cite[pages 431 ff.]{grecokiyek}).
\end{Defi}

\begin{Defi}
A matrix
\begin{equation}
\mathrm{HN}:=\left(%
\begin{array}{c}
  p_i \\
  c_i \\
  a_i \\
\end{array}%
\right)_{1 \le i < \infty} \nonumber
\end{equation}
with $p_i,c_i \in \nn \cup \{\infty\}$, $a_i \in k_i \setminus \{
0\}$ for every $i \in \nn$, is called an (abstract)
Hamburger-Noether tableau if it satisfies the following
properties:
\begin{itemize}
\item if $p_i = \infty$ for some $i \in \mathbb{N}$, then we have
$p_j=\infty$ and $c_j=c_i < \infty$ for every $j \ge i$;
    \item if $c_i=\infty$ for some $i \in \mathbb{N}$, then $c_j=\infty$ and $p_j=p_1<\infty$ for $j \in \mathbb{N}$.
    \item Assume that $c_1<\infty$. Then
    $c_{i+1}=\mathrm{gcd}(c_i,p_i)=1$ for every $i \in \nn$.
    \item We have $a_i=0$ if and only if $p_i=\infty$ or
    $c_i=\infty$, for every $i \in \mathbb{N}$.
\end{itemize}
If $p_1 < \infty$ and $c_1 < \infty$, then the tableau is said to
be non-degenerated.
\end{Defi}

\Nr Note that the Hamburger-Noether tableau of $(x,y)$ satisfies
the properties of the previous definition.

\Nr \label{614} Assume $\mathrm{HN}$ to be a non-degenerated
Hamburger-Noether tableau.
\begin{enumerate}
    \item An integer $i \in \{ 1, \ldots , l\}$ is called
    \emph{characteristic index} of $\mathrm{HN}$ if $i=1$ or if
    $c_{i+1}< c_i$.
    Let $1=i_1 < i_2 < \ldots < i_h$, $h:=h(\mathrm{HN}) \in
    \mathbb{N}$ be the characteristic indices of $\mathrm{HN}$. It is clear
    that $c_j=1$ for every $j \ge i_h +1$.
    \item Let us define $q_1 :=  p_1$, $q_j :=  p_{i_{j-1}+1}+ \ldots + p_{i_{j}}$ for
every $j \in \{2, \ldots , h \}$, $d_j :=  c_{i_j}$ for every $j
\in \{1, \ldots , h \}$. The sequence $Ch(\mathrm{HN}):=(d_1;q_1,
\ldots , q_h)$ is called the \emph{characteristic sequence} of
$\mathrm{HN}$. Notice that 
\[
d_i =  \mathrm{gcd} (d_{i-1},q_{i-1})=
\mathrm{gcd} (d_1, q_1, \ldots , q_{i-1})
\]
 for all $i \in \{2,
\ldots h \}$ and $d_{h+1} =  \mathrm{gcd} (d_h,q_h) = \mathrm{gcd}
(d_1, q_1, \ldots , q_h) = 1$. The sequence $d(\mathrm{HN}):=(d_1,
\ldots , d_{h+1})$ is called the \emph{divisor sequence} of
$\mathrm{HN}$. Notice also that
\begin{itemize}
    \item if $h=1$, then $d_1=c_1=1$ and $d_2=1$;
    \item if $h \ge 2$, then $d_1 \ge d_2 > d_3 > \ldots > d_h > d_{h+1}$
    and if $d_1=d_2$ then either $d_1 \mid q_1$ or $q_1 \mid d_1$.
\end{itemize}
    \item Furthermore, we set
    \[
n_i:=\frac{d_i}{d_{i+1}} ~ ~ ~ ~ \mathrm{~for~every~} i \in \{1,
\ldots , h \};
    \]
    the sequence $n(\mathrm{HN})=(n_1, \ldots,n_h)$ is called the
    $n$-sequence of $\mathrm{HN}$. We also set $r_0 := d_1$ and
\[
r_i :=  \sum_{j=1}^{i} q_j \frac{d_j}{d_i}
\]
for every $i \in \{1, \ldots , h \}$. The sequence
$r(\mathrm{HN}):=(r_0, \ldots , r_h)$ is called the
    \emph{semigroup sequence} of $\mathrm{HN}$.
\end{enumerate}

    \Nr If $\mathrm{HN}$ is a degenerated Hamburger-Noether tableau, then
    we will define $h(\mathrm{HN})=0$ and
    $Ch(\mathrm{HN})=n(\mathrm{HN})=r(\mathrm{HN})=(1)$.

\section{Curvettes} \label{sec:cur}

Let $f \in \mm_R$. Following \cite[page 397]{grecokiyek}, every $S
\supset R $ with $(fR)^S \ne S$ is said to be a locus point of
$fR$ (or of $f$); the set $\mathbb{L}(f)$ of locus points of $f$
is called the \emph{point locus} of $f$. Notice that
$\mathbb{L}(f)$ is an infinite set.

\begin{Defi}
Let $f \in \mm_R$. Set $S \supset R$ so that $S \notin
\mathbb{L}(f)$ and $S$ not be the intersection of two components
of the exceptional divisor. A \emph{curvette} at $S$ is defined to
be a normal-crossing curve $g \in S$ such that $(gS \cap R)^S$ is
a curve with no singularities at $S$ and not passing through any
other point $S^{\prime} \in \Omega (R)$ with $S \ne S^{\prime}$.
\end{Defi}

\begin{Prop} \label{prop:52}
Let $g \in S$ a normal crossing curve, where $S$ satisfies the
conditions of the above definition. Then there exists $h \in R$
irreducible with $gS \cap R = hR$ and $(hR)^S$ is normal-crossing
at $S$ with $(hR)^{S^{\prime}}=S^{\prime}$ for every $S^{\prime}
\in \Omega (R)$ with $S^{\prime} \nsupseteq S$ and $S^{\prime} \nsubseteq
S$.
\end{Prop}

\dem~ First assume that $S \in N_1(R)$, then there exists $p \in
\mathbb{P}_R$ generated by an irreducible homogeneous polynomial
$\overline{h} \in k_R[\overline{x},\overline{y}]$ (cf. \ref{2:2}).
Choose $h \in \mm_R^l$ with $\overline{h}=h \mathrm{~mod~}
\mm_R^{l+1}$. Without loss of generality, we assume that
$\overline{x}$ does not divide $\overline{h}$. Then the
exceptional divisor has the equation $xS$ and the strict transform
of $h$ in $S$ is $\frac{h}{x^l}S = (hR)^S$. Thus $\left (
x,\frac{h}{x^l} \right )$ is a regular system of parameters of
$S$. Inductively, it is easy to check this statement for every $S
\in N(R)$. Assume $xS$ is the equation of the exceptional divisor
in $S$. We have also that $g S \cap R \ne (0)$ and $g S \cap R \ne
\mm_R$: on the contrary, we would have $x \in g S$, which is a
contradiction. Hence $g S \cap R$ is a prime ideal of $R$
different from $0$ of height $1$. Since $(g S \cap R) S \cap R = g
S \cap R$, the transform $(g S \cap R)^S$ is a principal prime
ideal of $S$ with $(g S \cap R)^S \cap R = g S \cap R$ (by Lemma
\ref{lem:stricttransfisirred} (i)), and therefore $(g S \cap
R)^S=g S$. Moreover, for any other subring $S^{\prime}$ such that
$S^{\prime}$ is not infinitely near to $S$ and $S$ is not
infinitely near to $S^{\prime}$, again Lemma
\ref{lem:stricttransfisirred} shows us that $(g S \cap
R)^{S^{\prime}} = S^{\prime}$. \qed

\begin{Rem} \label{rem:53}
Notice that by Lemma \ref{lem:stricttransfisirred} (i) and
\cite[Chapter VII, (1.1)]{kiyek}, the strict transform $(gR)^S$ is
irreducible in $S$ and $gS \cap R=gR$.
\end{Rem}

The first result relating the intersection multiplicity of two
curves and their strict transforms---already defined in
\ref{section:intmult}---is the following (see \cite[Chapter VII,
(8.8)-(8.9)]{kiyek}):

\begin{Lem}[Intersection formula] \label{Lem:intfor}
Let $R \in \Omega (\mathcal{K})$ and $\{f,g \}$ be a regular
sequence in $R$. Then we have
\[
\iota_R (fR,gR) = \sum_{S \in N(R)} [S:R] \mathrm{ord}_S ((fR)^S)
\mathrm{ord}_S ((gR)^S).
\]
Moreover, we have
\[
\mathrm{ord}_R(f)=\sum_{S \in N_1 (R)} [S:R] \iota_S
((fR)^S,\mathfrak{m}_R S).
\]
\end{Lem}

Let us take now the proximity matrix $P_{\mathcal{C}}$ with
respect to the cluster $\mathcal{C}$ associated with the
resolution of $f$, and its inverse matrix
$Q_{\mathcal{C}}:=P_{\mathcal{C}}^{-1}$. We give now an
interpretation of the entries $q_{R,S}$, for $R,S \in
\mathcal{C}$, of the matrix $Q_{\mathcal{C}}$ in terms of
curvettes.

\begin{Theo} \label{prop:curvette}
Let $\mathcal{C}$ be the cluster of a resolution of a curve $f \in
R$.
\begin{itemize}
    \item[(i)] For every $S \in \mathcal{C}$, we have
    \[
\mathrm{ord}_S ((fR)^{S})= \sum_{T \succ S} [T:S] \mathrm{ord}_{T}
((fR)^{T}).
    \]
\end{itemize}
    \item[(ii)] Furthermore, Let $S \in \mathcal{C}$. For any $T \in N(R)$ and any curvette $g \in T$, we have
    \[
\mathrm{ord}_S ((gT \cap R)^{S})= \sum_{T^{\prime} \succ S}
[T^{\prime}:S] \mathrm{ord}_{T^{\prime}} ((gT \cap
R)^{T^{\prime}}).
    \]
\end{Theo}

\dem~Statement (ii) follows easily from (i), and this is a
consequence of \ref{Lem:intfor}. \qed

\begin{Cor} \label{cor:curvette}
Let $R^{\prime} \in N_n(R)$. Let $Q=Q_{\mathcal{C}}=(q_{S,T})$,
$S,T \in \mathcal{C}$ be the inverse of the proximity matrix
$P_{\mathcal{C}}$. The following statements hold:
\begin{itemize}
    \item[(i)]
     For any $S \in \mathcal{C}$, we have
    \[
q_{S,R^{\prime}} = \mathrm{ord}_S ((fR)^S).
    \]
    \item[(ii)] For any curvette $g \in T$, $T \in N(R)$ and any $S \in
    \mathcal{C}$, we have
    \[
q_{S,T}= \mathrm{ord}_S ((gT \cap R)^S).
    \]
\end{itemize}
\end{Cor}

\dem~ First of all, we reformulate the equation (i) in Theorem
\ref{prop:curvette} to have
\[
\sum_{T^{\prime}} [T:T^{\prime}] p_{T,T^{\prime}}
\mathrm{ord}_{T^{\prime}} ((fR)^{T^{\prime}}) =
\delta_{T,R^{\prime}}
\]
for all $T \in \mathcal{C}$, where
\[
\delta_{T,R^{\prime}} = \left \{
\begin{array}{ll}
    [T:R^{\prime}], & \mathrm{~if~} T=R^{\prime}; \\
    0, & \mathrm{~otherwise}. \\
\end{array}%
\right.
\]
We take now multiplication by $q_{S,T}$ for $S \in \mathcal{C}$
and sum over all $T \in \mathcal{C}$:
\[
\sum_T \sum_{T^{\prime}} [T:T^{\prime}] q_{S,T} p_{T,T^{\prime}}
\mathrm{ord}_{T^{\prime}} ((fR)^{T^{\prime}}) = \sum_{T} q_{S,T}
\delta_{T,R^{\prime}}.
\]
Since $Q$ is the inverse matrix of $P$, all terms cancel except for
those containing $R^{\prime}$, and we get
\[
q_{S,R^{\prime}}=\mathrm{ord}_{S} (fR)^S.
\]
The same argument works to prove the statement (ii) replacing
$R^{\prime}$ (resp. $fR$) by $T$ (resp. $gT \cap R$). \qed

\begin{Prop} \label{prop:curvette2}
Let be the matrix $M_{\mathcal{C}}:=Q_{\mathcal{C}}^{t} \cdot
\Delta_{\mathcal{C}}^{-1} \cdot Q_{\mathcal{C}}$. Let $T_1,T_2$ be
two points of $\mathcal{C}$. Then the $(T_1,T_2)$-entry of the
matrix $M_{\mathcal{C}}$ is equal to the intersection number
\[
\iota_R (g_1T_1 \cap R, g_2T_2 \cap R)
\]
of two curvettes $g_1,g_2$ of $T_1$ and $T_2$, respectively.
\end{Prop}

\dem~ It is just to consider the equalities
$-N_{\mathcal{C}}^{-1}=(P_{\mathcal{C}} \cdot \Delta_{\mathcal{C}}
\cdot P_{\mathcal{C}}^{t})^{-1}=Q_{\mathcal{C}}^{t} \cdot
\Delta_{\mathcal{C}}^{-1} \cdot Q_{\mathcal{C}}$, and the
intersection's formula \ref{Lem:intfor} applied to the cases
$fR=g_1T_1 \cap R$ and $gR = g_2 T_2 \cap R$. \qed

\section{Curvettes and approximations} \label{sec:approx}

    \Nr \label{notations:hensel} Let $R$ be complete. Let $\{x,y\}$ be a regular system of parameters of $R$. Let $f
    \in R$ be an analytically irreducible curve. The ring
    $S:=R/fR$ is an analytically irreducible local domain of
    dimension $1$ whose integral closure is a discrete valuation
    ring, which is a finitely generated $S$-module (cf. \cite[Chapter II, (3.17)]{kiyek}). Let us assume
    that $S$ does contain a perfect field $\mathbb{F}$. By
    Hensel's lemma, there exists a finite extension
    $\mathbb{F}^{\prime}$ of $\mathbb{F}$ such that $\mathbb{F}^{\prime} \subseteq
    \widehat{S}$ is a coefficient field in the completion
    $\widehat{S}$ of $S$ with respect to the Jacobson radical.
    Notice that $\mathbb{F}^{\prime}$ is uniquely determined: it
    is nothing but the integral closure of $\mathbb{F}$ in
    $\widehat{S}$. Since $\widehat{S} \cong \overline{\widehat{S}} \cong
    \widehat{\overline{S}}$ (where $\overline{\cdot}$ denotes
    ``integral closure" respect to the quotient field), we have $\mathbb{F}^{\prime} \subseteq
    \widehat{\overline{S}}$ and again by Hensel's lemma there
    exists a finite extension $\mathbb{F}^{\prime \prime}$ of
    $\mathbb{F}^{\prime}$ (uniquely determined as well) which is a
    coefficient field for $\widehat{\overline{S}}$. Consequently,
    for every uniformising parameter $t$ of $\overline{S}$ we have
    $\widehat{S} \to \mathbb{F}^{\prime \prime}[\![t]\!] \cong
    \widehat{\overline{S}}$, and so a natural morphism $\chi: R \to \mathbb{F}^{\prime
    \prime}[\![t]\!]$.

    \begin{Defi} \label{Defi:HNf}
The Hamburger-Noether tableau associated with an analytically
irreducible curve $f \in R$ is defined to be
\[
\mathrm{HN}(f;x,y):=\mathrm{HN}(\chi (x),\chi (y)).
\]
    \end{Defi}

\Nr Assume $f \in R$ to be residually rational, i.e., so that
$k_R$ is isomorphic to $\mathbb{F}$ and
$\mathbb{F}=\mathbb{F}^{\prime}=\mathbb{F}^{\prime \prime}$. Let
$w$ be discrete valuation given by the order function of $k_R
[\![t]\!]$. Set $\overline{x}:=\chi (x)$, $\overline{y}:= \chi
(y)$. If $\ox = \overline{0}$ (resp. $\yy=\overline{0}$), then an
easy reasoning shows that $f=ux$ (resp. $f=u^{\prime}y$), for
$u,u^{\prime}$ units in $R$. Let us assume that both $\ox$ and
$\yy$ are non-zero. Consider the Hamburger-Noether tableau
$\mathrm{HN}(x,y;f)$ of Definition \ref{Defi:HNf}. We have
$\ox=\omega_x t^{c_1}$ and $\yy = \omega_y t^{p_1}$ for
$\omega_x,\omega_y \in k_R \setminus \{0\}$. Note that if
$c_1=p_1$, then we write
$\omega^{\prime}=\frac{\omega_y}{\omega_x}$ and
$\mathrm{ord}(f)=c_1=p_1$; we set
$\yy^{\prime}:=\yy-\omega^{\prime} \ox$ and
$y^{\prime}:=y-\omega^{\prime} x$. Then $w(\ox) = c_1 <
w(\yy^{\prime})$ and there exists $\lambda \in k_R \setminus
\{0\}$ with $\mathrm{In}(f)=\lambda (y-\omega^{\prime}x)^{c_1}$.
If $c_1<p_1$, then $\mathrm{ord}(f)=w(\ox)=c_1$ and there exists
$\theta \in k_R \setminus \{0\}$ with $\mathrm{In}(f)=\theta
y^{c_1}$. Also, if $c_1>p_1$, then $\mathrm{ord}(f)=w(\yy)=p_1$
and there is $\theta^{\prime} \in k_R \setminus \{0\}$ with
$\mathrm{In}(f)=\theta^{\prime} x^{p_1}$. Hence we assume in every
case that $\mathrm{In}(f)=(\lambda x + \mu y)^{\min (c_1,p_1)}$,
with $\lambda, \mu \in R$ not vanishing simultaneously (just by
multiplying with an element of $k_R \setminus \{0\}$). If
$\lambda=0$, then $f$ is said to be $y$-regular; if $\mu=0$, then
$f$ is said to be $x$-regular.
\medskip

Let $f \in R$ be an analytically irreducible residually rational
$y$-regular curve with $f \ne u y$ for some unit $u \in R$. Set
$\mathrm{HN}:=\mathrm{HN} (f;x,y)$, $h:=h(\mathrm{HN}(f;x,y))$,
$r=r(\mathrm{HN}(f;x,y))$, $d=d(\mathrm{HN}(f;x,y))$. We adapt
some results proven for algebroid curves in \cite{russell} to our
more general case. Next lemma corresponds to \cite[Lemma
2.10]{russell}.

\begin{Lem} \label{Lem:616}
Let be the two following Hamburger-Noether tableaux
\begin{equation}
\mathrm{HN}=\left(%
\begin{array}{c}
  p_i \\
  c_i \\
  a_i \\
\end{array}%
\right)_{1 \le i <\infty}, ~ ~ ~ ~
\mathrm{HN}^{\prime}=\left(%
\begin{array}{c}
  p^{\prime}_i \\
  c^{\prime}_i \\
  a^{\prime}_i \\
\end{array}%
\right)_{1 \le i <\infty}, \nonumber
\end{equation}
and let $s \in \mathbb{N}$. The following statements are
equivalent:
\begin{enumerate}
    \item We have $\frac{p_j}{c_1} =
    \frac{p^{\prime}_j}{c^{\prime}_1}$ for every $j \in \{1, \ldots , s
    \}$.
    \item We have $\frac{p_j}{c_i} =
    \frac{p^{\prime}_j}{c^{\prime}_i}$ for all $i,j \in \{1, \ldots , s
    \}$.
    \item We have $\frac{p_j}{c_j} =
    \frac{p^{\prime}_j}{c^{\prime}_j}$ for every $j \in \{1, \ldots , s
    \}$.
\end{enumerate}
Moreover, each of these conditions implies that
\[
\frac{p_j}{c_{s+1}}=\frac{p^{\prime}_j}{c^{\prime}_{s+1}}  ~ ~ ~ ~
\mathrm{~and~} ~ ~ ~ ~ \frac{c_j}{c_{s+1}} =
\frac{c^{\prime}_j}{c_{s+1}} ~ ~ ~ \mathrm{~for~every~} j \in \{
1, \ldots ,s \}.
\]
\end{Lem}

\begin{Defi}
Let be the two Hamburger-Noether tableaux
\begin{equation}
\mathrm{HN}=\left(%
\begin{array}{c}
  p_i \\
  c_i \\
  a_i \\
\end{array}%
\right)_{1 \le i <\infty}, ~ ~ ~ ~
\mathrm{HN}^{\prime}=\left(%
\begin{array}{c}
  p^{\prime}_i \\
  c^{\prime}_i \\
  a^{\prime}_i \\
\end{array}%
\right)_{1 \le i <\infty}. \nonumber
\end{equation}
 We set
\begin{align*}
S(\mathrm{HN},\mathrm{HN}^{\prime}) & =   \{0 \} \cup \Big \{ j \in \mathbb{N} \mid \frac{p_i}{c_i} = \frac{p^{\prime}_i}{c^{\prime}_i}, ~ ~ a_i=a^{\prime}_i ~ ~ \mathrm{~for~} i \le j \Big \}.\\
s(\mathrm{HN},\mathrm{HN}^{\prime}) & =  \mathrm{sup}
(S(\mathrm{HN},\mathrm{HN}^{\prime})).
\end{align*}
Notice that if $\mathrm{HN}=\mathrm{HN}^{\prime}$, then we have
$s(\mathrm{HN},\mathrm{HN}^{\prime})=\infty$.
\end{Defi}

\begin{Lem} \label{Lem:621}
Let $f,g \in R$ be two analytically irreducible and residually
rational curves. Let
\begin{equation}
\mathrm{HN}:=\mathrm{HN}(f;x,y)=\left(%
\begin{array}{c}
  p_i \\
  c_i \\
  a_i \\
\end{array}%
\right)_{1 \le i <\infty}, ~ ~ ~ ~
\mathrm{HN}^{\prime}:=\mathrm{HN}^{\prime}(g;x,y)=\left(%
\begin{array}{c}
  p^{\prime}_i \\
  c^{\prime}_i \\
  a^{\prime}_i \\
\end{array}%
\right)_{1 \le i <\infty} \nonumber
\end{equation}
be the Hamburger-Noether tableaux of $f$, resp. of $g$. Set
$s:=s(\mathrm{HN},\mathrm{HN}^{\prime})$. Then we have
\begin{align*}
\iota_R (fR,gR) & = \sum_{i=1}^{s} p_i c^{\prime}_i + \mathrm{min} (\{ p_{s+1} c^{\prime}_{s+1},p^{\prime}_{s+1} c_{s+1} \})  \\
 & =  \sum_{i=1}^{s} p^{\prime}_i c_i + \mathrm{min} (\{ p_{s+1} c^{\prime}_{s+1},p^{\prime}_{s+1} c_{s+1} \}).
\end{align*}
\end{Lem}

\dem~The reasoning is much more similar as that for algebroid
curves in \cite[Theorem 3.3]{russell}. \qed

\begin{Defi}
Let $\mu \in \mathbb{N}$. A Hamburger-Noether tableau
$\mathrm{HN}^{\prime}$ is called a $\mu$-th approximation to
$\mathrm{HN}$ if
\begin{enumerate}
    \item $s(\mathrm{HN},\mathrm{HN}^{\prime})=\mu-1$;
    \item $c^{\prime}_{\mu}=1$;
    \item $p^{\prime}_{\mu} c_{\mu} \ge p_{\mu}$.
\end{enumerate}
\end{Defi}

\begin{Rem} \label{Rem:pc}
Let $\mathrm{HN}^{\prime}$ be a $\mu$-approximation to
$\mathrm{HN}$. In this case we have
\[
p^{\prime}_i = \frac{p_i}{c_{\mu}} \ \ \mathrm{~and~} \ \
c^{\prime}_i = \frac{c_i}{c_{\mu}}
\]
for every $i \in \{1, \ldots , \mu -1\}$.
\end{Rem}

\begin{Defi}
A curve $g \in R$ is called a $\mu$-approximation to $f$ if $g$ is
analytically irreducible, residually rational and if
$\mathrm{HN}(g;x,y)$ is a $\mu$-th approximation to
$\mathrm{HN}(f;x,y)$.
\end{Defi}

\begin{Prop} \label{prop:709}
Let $g \in R$ be a $\mu$-th approximation to $f$. There exists a
curvette $h \in T$ for some $T \in N(R)$ with $hT \cap R=gR$.
Conversely, given a curvette $h \in T$ for some $T \in N(R)$,
there exists a $\mu$-th approximation $g$ to $f$ such that $hT
\cap R = gR$.
\end{Prop}

\dem~ Let $g \in R$ be a $\mu$-th approximation to $f$. Assuming
$\mathrm{ord}_R(g)=l$ and $\yy$ does not divide $h \mathrm{~mod~}
\mm_R^{l+1}$, then the strict transform of $g$ in $T$ is
$\frac{g}{y^l}T=(gR)^T$, and $\left ( y,\frac{g}{y^l} \right )$ is
a regular system of parameters of $T$. Then $h:=\frac{g}{y^l}$ is
a curvette in $T$ with $hT \cap R = g$ (by the same reasoning as
in the proof of Proposition \ref{prop:52}). Conversely, if $h \in
T$ is a curvette for some $T \in N(R)$, again by Proposition
\ref{prop:52} there exists an irreducible element $g \in R$ with
$hT \cap R=gR$. To prove that $g \in R$ is a $\mu$-th
approximation to $f$, since $g$ is analytically irreducible and
residually rational (see Remark \ref{rem:53}), it suffices to use
lemmas \ref{Lem:616} and \ref{Lem:621} and to argue like in
\cite[pages 59--60]{russell}.
 \qed

\begin{Rem}
Let $g \in R$ be a $\mu$-th approximation to $f$. We have
\[
\iota_R (fR,gR)=\sum_{i=1}^{\mu} p_i \frac{c_i}{c_{\mu}}.
\]
\end{Rem}

\begin{Theo} \label{prop:711}
Let $f \in R$ be an analytically irreducible residually rational
$y$-regular curve with $f \ne u y$ for some unit $u \in R$. Let
$\mathcal{C}_s$ be the cluster associated with the minimal
resolution $\pi=\pi_1 \circ \ldots \circ \pi_s$ of $f$. Let $T_j$
be a non-singular point of the $j$-th component of the exceptional
divisor of $\pi_s$, for $1 \le j \le s$. Let $g^{(j)} \in T_j$ be
an analytically irreducible and residually rational curve. The
following assertions are equivalent:
\begin{enumerate}
    \item $g^{(j)}$ is a curvette on $T$;
    \item $\iota_R (fR,g^{(j)}T \cap R)=r_j$.
\end{enumerate}
\end{Theo}

\dem~By Proposition \ref{prop:709}, the curve $g^{(j)}$ is a
$\mu$-th approximation to $f$, where $\mu$ is the $j$-th
characteristic index of $\mathrm{HN}(f;x,y)$, i.e., we have $\mu
=i_j$, for $j \in \{1, \ldots ,h\}$. We consider two cases: (i)
Assume that $j=1$ and $d_1 \mid r_1$. The proof of $(2)
\Rightarrow (1)$ follows from Lemma \ref{Lem:621}, and the
converse is also easy. (ii) Assume that $j > 1$, or $j=1$ and $d_1
\nmid r_1$. Since $\mu =i_j$, by \ref{614} we have
\[
\sum_{i=1}^{\mu} p_i \frac{c_1}{c_{\mu}} = \sum_{i=1}^{j} q_i
\frac{d_i}{d_j} =r_j
\]
and $d_j=c_{\mu}$. By \ref{614} and Remark \ref{Rem:pc}, we see
that (2) follows from (1). Conversely, let us assume that (2)
holds. Let
$$
\mathrm{HN}^{\prime}:=\mathrm{HN}^{\prime}(g^{(j)};x,y)=\left(%
\begin{array}{c}
  p^{\prime}_i \\
  c^{\prime}_i \\
  a^{\prime}_i \\
\end{array}%
\right)_{1 \le i <\infty} .\nonumber
$$
We set now $s=s(\mathrm{HN},\mathrm{HN}^{\prime})$, and we proceed
as in the proof of Proposition \ref{prop:709}. \qed

\section*{Acknowledgements}

The author wishes to express his gratitude to Prof. Dr. Karlheinz
Kiyek for teaching him carefully the basics of the theory of
two-dimensional regular local rings. He is also thankful to the
referee of this paper for several helpful comments.

\end{document}